\documentclass[12pt,epsfig]{amsart}

\usepackage{amssymb, adjustbox, enumerate, amsbsy, stmaryrd}
\usepackage[dvipsnames]{xcolor}
\usepackage{amsmath,wasysym}
\usepackage{comment}
\usepackage{caption}
\usepackage{stackrel}
\usepackage[mathscr]{eucal}
\usepackage{amsthm}
\usepackage{graphicx}
\usepackage {amscd}
\usepackage {epic}
\usepackage[all,color]{xy}
\usepackage[alphabetic]{amsrefs}
\usepackage{extarrows}  % \xlseftrightarrow
\usepackage{color}
\usepackage{braket}

\usepackage{graphpap, color}

\usepackage[mathscr]{eucal}
\usepackage{mathrsfs}%,mathabx}
\usepackage{graphicx}

\usepackage{tikz}
\usepackage{tikz-cd}
\usepackage{tikz-3dplot}
\newtheorem{theorem}{Theorem}

\newtheorem{claim}[theorem]{Claim}
\newtheorem{proposition}[theorem]{Proposition}
\newtheorem{lemma}[theorem]{Lemma}
\newtheorem{definition}[theorem]{Definition}
\newtheorem{corollary}[theorem]{Corollary}
\newtheorem{remark}[theorem]{Remark}
\newtheorem{example}[theorem]{Example}
\newtheorem{question}[theorem]{Question}

\numberwithin{equation}{section} 
\numberwithin{theorem}{section}

\def\cF{\mathcal{F}}

\def\cT{\mathcal{T}}

\def\CC{\mathbb{C}}

\def\GG{\mathbb{G}}

\def\bbN{\mathbb{N}}
\def\PP{\mathbb{P}}

\def\bbR{\mathbb{R}}
\def\bbZ{\mathbb{Z}}

\newcommand{\restrict}[2]{\ensuremath{\left. #1 \right|_{#2}}}

\def\conv{\mathrm{Conv}}
\def\dd{\mathrm{partial}}

\def\Hom{\mathrm{Hom}}

\def\SL{\mathrm{SL}}

\def\supp{\mathrm{Supp}}

\def\vol{\mathrm{vol}}

\def\dd{\partial}

\def\bs{\boldsymbol}

\def\vep{\varepsilon}

\newcommand{\reduce}[1]{\scalebox{1}{\ensuremath{#1}}}

\newcommand{\norm}[1]{\ensuremath{\left| #1 \right|}}
\newcommand{\Norm}[1]{\ensuremath{\left\| #1 \right\|}}

\DeclareMathOperator*{\smalloplus}{\reduce{\bigoplus}}

\DeclareMathOperator*{\smallcap}{\reduce{\bigcap}}
\DeclareMathOperator*{\smallcup}{\reduce{\bigcup}}

\input{xy}
\xyoption{all}

\title[Chow stability of uniform K-stable toric varieties]{Asymptotic Chow stability of uniformly \\K-stable toric varieties}
\author{King leung Lee}
\address{IMAG, 
Universit\'{e} de Montpellier, 499-554 Rue du Truel 9, 34090 Montpellier, France}
\email{king-leung.lee@umontpellier.fr}

\author{Naoto Yotsutani }
\address{Department of Mathematics, Faculty of Science, Shizuoka University, $836$ Ohya, Suruga-ku, Shizuoka-shi, Shizuoka, $422$-$8529$, Japan}
\email{yotsutani.naoto@shizuoka.ac.jp}

\address{IMAG, Universit\'{e} de Montpellier, $499$-$554$ Rue du Truel $9$, $34090$ Montpellier, France}
\email{naoto.yotsutani@umontpellier.fr}

\makeatletter
\@namedef{subjclassname@2020}{%
  \textup{2020} Mathematics Subject Classification}
\makeatother

\subjclass[2020]{Primary 51M20; Secondary  53C55, 14M25}

\keywords{Triangulations, the Bramble-Hilbert Lemma, toric varieties, convex geometry, }

\date{\today}

\begin{document}

\maketitle

\begin{abstract}
For a polarized toric variety, we provide a sufficient criterion ensuring that a uniformly K-stable polarized toric variety $(X,L)$ is asymptotically Chow polystable, under the assumption that the obstruction to asymptotic Chow semistability (the Futaki-Ono invariant) vanishes. Our approach is based on a detailed study of triangulations of neighborhoods of the vertices of the associated moment polytope $\Delta$. 
%when the obstruction of asymptotic Chow semistability (the Futaki-Ono invariant) vanishes by studying the triangulation of the neighborhoods of the vertices of the corresponding moment polytope $\Delta$. 
As an application, we prove that every uniformly K-stable polarized smooth toric variety $(X,L)$ with vanishing Futaki-Ono invariant is asymptotically Chow polystable.
\end{abstract}

\maketitle

%Sec1
\section{Introduction}
One of the central problems in K\"ahler geometry is the existence of constant scalar curvature K\"ahler (cscK) metrics in the first Chern class of a polarized variety $(X,L)$.
The Yau-Tian-Donaldson conjecture asserts that the existence of a cscK metric in $c_1(L)$ is equivalent to the uniform K-stability of $(X,L)$.
In particular, this conjecture has been established by Delcroix \cite{Del23} and Li \cite{Li22} for polarized {\emph{spherical varieties}}, which can be viewed as a natural generalization of polarized toric varieties.

In this paper, we study uniform K-stability, which has been extensively investigated in \cite{BHJ17, His20}. We emphasize that the main idea underlying our approach to polarized toric varieties is closely related to the following inequality, which originally appeared in \cite[Proposition $5.1.2$]{Don02}:
\[
\mathcal L_a(f)\geq \lambda \int_{\partial \Delta}f(\bs x)d\sigma,
\]
where 
\[
\mathcal L_a(f):=\int_{\partial \Delta}f(\bs x)d\sigma-a\int_\Delta f(\bs x)d\mu, \qquad\text{with} \quad a=\frac{\vol(\partial \Delta, d\sigma)}{\vol(\Delta)}.
\]
%with $a=\vol(\partial \Delta, d\sigma)/\vol(\delta)$. 
When $\lambda=0$, this inequality reduces to the condition of K-semistability for a polarized toric variety.
%The aim of this paper is to introduce the notion of {\emph{$\lambda$-stability}} for polarized toric varieties which was originally appeared in \cite{Don02}, Proposition $5.1.2$. Then, we see that $\lambda$-stability is a natural generalization of {\emph{uniform K-stability}} which was well stablished in \cite{BHJ17, His20}. 
%In particular, $\lambda$-stability coincides with K-semistability when a real number $\lambda$ equals to $0$. See, Definition \ref{def:Lambda_Sta} for the definition of $\lambda$-stability. 

Another important notion of GIT-stability in K\"ahler geometry is Chow stability (see Section \ref{sec:Chow} for the definition of asymptotic Chow stability).
Let $(X,L)$ be an $n$-dimensional polarized manifold, and let $\mathrm{Aut}(L)$ denote the group of bundle automorphisms of the ample line bundle $L$.
Since $\mathrm{Aut}(L)$ naturally contains $\CC^\times:=\CC \setminus \set{0}$ as a subgroup acting by fiberwise scalar multiplication, we define
%which acts as fiber multiplications, we set 
\[
\mathrm{Aut}(X,L):=\mathrm{Aut}(L)/\CC^\times.
\]
Then any element of $\mathrm{Aut}(X,L)$ induces an automorphism of $X$, and hence $\mathrm{Aut}(X,L)$ may be regarded as a Lie subgroup of $\mathrm{Aut}(X)$. 

In \cite{Don01}, Donaldson proved that if $(X,L)$ admits a cscK metric in $c_1(L)$ and $\mathrm{Aut}(X,L)$ is discrete, 
then $(X,L)$ is asymptotically Chow stable.
This result was partially extended by Mabuchi to the case where $\mathrm{Aut}(X,L)$ is not discrete. 
More precisely, Mabuchi \cite{Mab05} proved that if $(X,L)$ admits a cscK metric in $c_1(L)$, then $(X,L)$ is 
asymptotically Chow polystable provided that the obstruction to asymptotic Chow semistablity vanishes.
%$(X,L)$ satisfies the vanishing condition for the obstruction to asymptotic Chow semistablity.
Subsequently, Futaki showed that Mabuchi's condition is equivalent to the vanishing of a family of integral invariants $\mathcal{F}_{\mathrm{Td}^1}, \dots, \mathcal{F}_{\mathrm{Td}^n}$ introduced in \cite{Fut04},
where $\mathrm{Td}^i$ denotes the $i$-th Todd polynomial. In particular, $\mathcal{F}_{\mathrm{Td}^1}$ coincides with the classical Futaki invariant up to a non-zero constant multiple. %multiplicative constant.

More precisely, recall that the classical Futaki invariant is the map $f:\mathfrak h(X)\to \CC$ defined as follows.
Let $(X,\omega)$ be a compact K\"ahler manifold with K\"ahler form $\omega$, and let $\mathfrak h(X)$ denote the Lie algebra of holomorphic vector fields on $X$.
Denote by $s(\omega)$ the scalar curvature of $\omega$, and let $(g^{i\bar{j}})$ %_{i\bar{j}}$ 
be the inverse of the K\"ahler metric $(g_{i\bar{j}})$. %_{i\bar{j}}$.
The complex Laplacian associated with $\omega$ is given by 
\[
\Delta_\omega:=-g^{i\bar{j}}\partial_i \bar{\partial}_j.
\]
Let $h_\omega$ be the real-valued function determined (up to an additive constant) by
\[
s(\omega)-%\left( \int_X s(\omega)\omega^n \biggl/ \int_X \omega^n\right)
\frac{\int_X s(\omega)\omega^n }{ \int_X \omega^n}
=-\Delta_\omega h_\omega.
\]
%up to addition of a constant. 

The Futaki invariant is then defined by
%$f:\mathfrak h(X) \to \CC$ is defined by
\[
f(v):=\int_X v (h_\omega) \omega^n,\qquad v\in \mathfrak h(X).
\]
It is well known that $f(v)$ provides an obstruction to the existence of a cscK metric in the K\"ahler class $[\omega]$.
%its first Chern class.

Let $\mathfrak h_0(X)$ denote the Lie subalgebra in $\mathfrak h(X)$ corresponding to $\mathrm{Aut}(X,L)$, 
%where $\mathfrak h_0(X)$ 
consisting of holomorphic vector fields with non-empty zero set. It is known that $\mathcal F_{\mathrm{Td}^1}$ coincides with the restriction $f|_{\mathfrak h_0(X)}$ up to multiplication by a non-zero constant.
Since the family of integral invariants $\mathcal{F}_{\mathrm{Td}^1}, \dots, \mathcal{F}_{\mathrm{Td}^n}$ extends the classical Futaki invariant, %are a generalization of the classical Futaki invariant, 
we refer to them as the {\emph{higher Futaki invariants}}. Combining Mabuchi's result \cite{Mab05} with Futaki's observation \cite{Fut04}, we obtain the following partial generalization of Donaldson's theorem \cite{Don01}.
%Thm
\begin{theorem}[Mabuchi-Futaki]
Let $(X,L)$ be an $n$-dimensional polarized manifold. Assume that the higher Futaki invariants $\mathcal{F}_{\mathrm{Td}^i}$ vanish for all $i=1,\dots, n$, and that $(X,L)$
admits a cscK metric in $c_1(L)$. Then $(X,L)$ is asymptotically Chow polystable.
\end{theorem}
In \cite{Ono11}, Ono reformulated the obstruction given by the higher Futaki invariants to 
%$\mathcal{F}_{\mathrm{Td}^i}$ to 
asymptotic Chow semistablity for polarized toric varieties %in terms of the associated moment polytope $\Delta$ 
as follows. Let $\Delta\subset M_\bbR\cong \bbR^n$ be an integral Delzant polytope, and let $E_\Delta(t)$ be its Ehrhart polynomial: %of $\Delta$ in the form of
\begin{equation*}
E_{\Delta}(t)=\vol(\Delta)t^n+\frac{\vol(\partial \Delta, \sigma)}{2}t^{n-1}+\sum_{j=0}^{n-2}E_{\Delta,j} t^j, %\quad \text{with} \quad  E_\Delta(k)=|k\Delta\cap \bbZ^n|
\end{equation*}
so that $E_\Delta(k)=|k\Delta\cap \bbZ^n|$ for every positive integer $k$. Similarly, there exists an $\bbR^n$-valued polynomial $\boldsymbol s_\Delta(t)$, called the {\emph{sum polynomial}} of $\Delta$, given by
\begin{align}\label{eq:SumPoly}
\boldsymbol s_\Delta(t) & =t^{n+1}\int_\Delta \bs x\, {d\mu}+\frac{t^{n}}{2}\int_{\partial \Delta} \bs x\, d\sigma+\sum_{j=1}^{n-1}t^j\bs s_{\Delta,j}, \quad \text{so that} \quad
\bs s_\Delta(k)=\sum_{\bs a\in k\Delta\cap \bbZ^n}\bs a.
\end{align}
Here $\bs x=(x_1, \dots, x_n)$. The measure $d\sigma$ denotes the $(n-1)$-dimensional Lebesgue measure on $\partial \Delta$, defined as follows. Let $h_i(\bs x)=\braket{\bs x,  \bs v_i}+a_i$ be the defining equation of a facet $F_i$ of $\Delta$, and let ${d\mu}=dx_1\wedge \dots \wedge dx_n$ be the standard volume form on $\bbR^n$. On each facet $F_i=\set{\bs x \in \Delta| h_i(\bs x)=0}\subset \partial \Delta$, the induced measure is determined by %$(n-1)$-dimensional Lebesgue measure $d\sigma_i=d\sigma|_{F_i}$ is defined by
\[
{d\mu}=\pm d\sigma_i\wedge dh_i.
\]
In \cite[Corollary $1.5$]{Ono11}, Ono introduced the invariants 
%$\mathcal F_{\Delta,j}$ which is the $\bbR^n$-valued polynomials
\begin{equation}\label{eq:Ono-inv}
\mathcal F_{\Delta,j}:=\vol(\Delta)\bs s_{\Delta,j}-E_{\Delta, j-1}\int_\Delta \bs x \, {d\mu}\in \bbR^n, \qquad j=1,\dots, n.
\end{equation}
He further proved that if the associated polarized toric manifold $(X_\Delta, L_\Delta)$ is asymptotically Chow semistable, then $\mathcal F_{\Delta, j}=0$ for all $j=1, \dots, n$. 
It was conjectured that the complex linear span of the $\mathcal F_{\Delta, j}$ coincides with that of the higher Futaki invariants
%the linear hull of $\mathcal F_{\mathrm{Td}^p}$ in $\CC^n$ 
(see, \cite[Conjecture $1.6$]{Ono11}), namely,
\begin{equation}\label{eq:LinHull}
\mathrm{Lin}_\CC\Set{\mathcal F_{\Delta,j}: ~ j=1, \dots , n}=\mathrm{Lin}_\CC\Set{\restrict{\mathcal F_{\mathrm{Td}^p}}{{\CC^n}}: ~ p=1, \dots , n} \subset \CC^n.
\end{equation}
%for a polarized toric manifold $(X_\Delta, L_\Delta)$. We remark that the equality in $\eqref{eq:LinHull}$ was justified by Futaki in \cite{Fut12}.
This conjecture was later proved by Futaki \cite{Fut12}.
By Proposition \ref{prop:FO_inv}, $\cF_{\Delta,j}=0$ for all $j=1, \dots, n$ is equivalent to 
\[
\frac{1}{|k\Delta\cap \bbZ^n|}\sum_{\bs p\in k\Delta\cap \bbZ^n}\ell(\bs p)-\frac{1}{\vol(k\Delta)}\int_{k\Delta}\ell(\bs x){d\mu}=0
\]
for every positive integer $k$ and every affine function $\ell(\bs x)$. %Therefore, we define the Futaki-Ono invariant to be 
This observation motivates the following definition of the {\emph{Futaki-Ono invariant}}:
\[
FO(\ell;k):=\frac{1}{|k\Delta\cap \bbZ^n|}\sum_{\bs p\in k\Delta\cap \bbZ^n}\ell(\bs p)-\frac{1}{\vol(k\Delta)}\int_{k\Delta}\ell(\bs x){d\mu}.
\]
%for a positive integer $k$ and an affine function $\ell(x)$.
Combining the existence theorem of cscK metrics for polarized $G$-varieties due to Delcroix \cite{Del23}, Chen-Cheng \cite{CC18, CC21}, and Hisamoto \cite{His20}, one obtains the following result.
%Thm
\begin{theorem}\label{thm:Delcroix}
Let $(X,L)$ be an $n$-dimensional uniformly K-stable polarized toric manifold with associated polytope $\Delta\subset M_\bbR\cong \bbR^n$.
If the Futaki-Ono invariant $FO(\ell;k)$ vanishes for every positive integer $k$ and every affine function $\ell(\bs x)$, then $(X,L)$ is asymptotically Chow polystable.
\end{theorem}
This naturally leads to the following question concerning the GIT-stability of polarized toric varieties, which was already raised in the second author's earlier work \cite{Yot18}.
% Problem
\begin{question}[cf. \cite{Yot18}, Section $2$]\label{prob:GIT}
Can one give a direct combinatorial proof of Theorem \ref{thm:Delcroix} using techniques from toric geometry?
\end{question}
In \cite{Yot18}, we outlined a strategy for proving Theorem \ref{thm:Delcroix} combinatorially, but were unable to complete 
the argument because of certain technical difficulties (see, \cite[Remark $3.4$]{Yot18}).  
One of the main goals of the present paper is to resolve Question \ref{prob:GIT}. 
We approach this problem by studying the vertex cones of the polytope $\Delta$. 

For a vertex $\bs p\in \mathcal V(\Delta)$, recall that the {\emph{vertex cone}} $C(\bs p)$ is defined by
\[
C(\bs p):=\set{\bs p+ \lambda (\bs x-\bs p) | \bs x \in \Delta, ~~ \lambda \in \mathbb R_{\geq 0}}.
\]
See \cite[$\S 11.3$]{BR15} for further details.
 
In Sections \ref{sec:SmallPolytope} and \ref{sec: sct}, we introduce the notions of {\emph{small polytopes}} and {\emph{semi-canonical traiangulations}} (see Definition \ref{def:SmallPolytope} and Definition \ref{def:SCT}, respectively). To this end, we prove the following main result.

%See, Sections \ref{sec:Ideal} and \ref{sec:SmallPolytope} for the definitions of an ideal triangulation and small polytope, respectively.
%Thm
\begin{theorem}[See, Theorem \ref{theo: main theorem 2}] \label{thm:main}
Let $(X,L)$ be a polarized toric variety with integral polytope $\Delta$. Suppose that: 
\begin{enumerate}
    \item $(X,L)$ is uniformly K-stable;
    \item all Futaki-Ono invariants vanish;
    \item every vertex cone of $\Delta$ admits a semi-canonical triangulation which is small. 
\end{enumerate}
Then $(X,L)$ is asymptotically Chow polystable.
\end{theorem}

Combining Theorem \ref{thm:main}  and Corollary \ref{standard cone is small}, we obtain the following resolution of Question \ref{prob:GIT}.
%this particular triangulation to an $n$-dimensional toric manifold, we solve Question \ref{prob:GIT} as follows.
%Cor
\begin{corollary}[See, Corollary \ref{coro: K stable manifold is Chow stable}]\label{cor:main}
Let $(X,L)$ be a uniformly K-stable polarized toric manifold with vanishing Futaki-Ono invariants.
Then, $(X,L)$ is asymptotically Chow polystable.
\end{corollary}
%Since uniform K-stability implies K-semistability, our result generalizes the known result (Theorem \ref{thm:Delcroix}) under an assumption weaker (i.e., K-semistability) than the original one (i.e., uniform K-stability). Hence, this is an advantage of our combinatorial approach. 
For the reader's convenience, we provide several explicit examples of asymptotically Chow polystable toric varieties in Section \ref{sec:Examples} by applying Corollary \ref{cor:main}. In particular, we construct a three-dimensional non-reflexive singular toric variety that is asymptotically Chow polystable. 

\vskip 7pt

\noindent{\bfseries Organization of the paper.}
%This paper is organized as follows: 
In Section \ref{sec:Prelim}, we review basic facts concerning toric varieties and Chow stability.
We also introduce the notion of ideal triangulations and recall criteria for uniform K-stability in the toric setting.
In Section \ref{sec:SmPolyUnifK}, we study small polytopes in greater detail and introduce the notion of Uniform K-stability of polarized toric varieties.
In Section \ref{sec:Criteria}, we define the Futaki-Ono invariants, which provide an obstruction to asymptotic Chow semistability for polarized toric varieties, and prove Theorem \ref{thm:main}.
In Section \ref{sec: sct}, we introduce semi-canonical triangulations and show that they are ideal; 
as a consequence, we prove Theorem \ref{thm:main} and Corollary \ref{cor:main}. In Section \ref{sec:SymWeak}, we construct a family of asymptotically Chow polystable toric varieties that are not necessarily Fano, illustrating an application of Theorem \ref{thm:main}. 
Section \ref{sec:Examples} contains several examples of small polytopes that highlight the combinatorial
and technical aspects of our theory. 
%After, we list several examples of small polytopes in Section \ref{sec:Examples} that would be helpful to illustrate combinatorial properties and technical features of thesepolytopes. 
%Finally, in the appendix, we provide a detailed proof of Lemma 7.1 from \cite{Lee25}, whose original proof is not fully explicit.

\vskip 7pt

\noindent{\bfseries Acknowledgements.}
This work was supported by the first author's ANR-$21$-CE$40$-$0011$ JCJC project MARGE, and by the second author's JSPS KAKENHI Grants JP$22$K$03316$ and JP$24$K$0252$, as well as NSFC Grant $12571058$.%(Grant-in-Aid for Scientific Research (C)). 

The first author would like to express his sincere gratitude to Professors Thibaut Delcroix, Yat Tin Chow, and Kwok Kun Kwong for many helpful discussions. He is also grateful to Virginie Iackle for her warm encouragement.
Finally, he would like to thank Professor Xiaowei Wang for carefully reading an earlier version of this manuscript and providing valuable comments and suggestions.
%The first author deeply appreciates Professors Thibaut Delcroix, Yat Tin Chow and Kwok Kun Kwong for their helpful discussions. Also, he would like to thank Virginie Iackle for her warm encouragement. Finally, the first author deeply appreciates Professor Xiaowei Wang for the revision of the work.

%Sec2
\section{Toric varieties and simplex triangulations}\label{sec:Prelim}
In this section, we briefly review geometric invariant theory and toric geometry.
We also introduce several special classes of triangulations of integral polytopes. 

%Sec2.1
\subsection{Chow stability and toric varieties}\label{sec:Chow}
%Firstly, we recall the notion of Chow stability. 
We begin by recalling the notion of Chow stability.
For further details, we refer the reader to \cite{Lee25, LLSW19, Ono11, Yot16}.

Let $G$ be a reductive algebraic group and $V$ a finite dimensional complex vector space. 
Suppose that $G$ acts linearly on $V$.
For a point $\bs v\in \PP(V)$, let $\bs v^{\ast}\in V$ be a representative, and denote by $\mathcal{O}_G(\bs v^{\ast})$ 
the $G$-orbit of $\bs v^\ast$ in $V$.
%Let us denote a point $v^{\ast}$ in $V$ which is a representative of $v\in \PP(V)$.
%Let $\mathcal{O}_G(v^{\ast})$ be the $G$-orbit in $V$. 
Then:
\begin{enumerate}
\item $\bs v^{\ast}$ is called \emph{$G$-semistable} if the Zariski closure $\overline{\mathcal O_G(\bs v^{\ast})}$ does not contain the origin, i.e.,  $0\notin \overline{\mathcal O_G(\bs v^{\ast})}$.
\item $\bs v^{\ast}$ is called \emph{$G$-polystable} if
$\mathcal{O}_G(\bs v^{\ast})$ is closed.
\end{enumerate}
Analogously, a point $\bs v\in \PP(V)$ is said to be $G$-polystable (resp. $G$-semistable) if any representative $\bs v^\ast$ is $G$-polystable (resp. $G$-semistable).
%Rem
\begin{remark}\rm
The closure of $\mathcal O_G(\bs v^{\ast})$ in the Euclidean topology coincides with the Zariski closure $\overline{\mathcal O_G(\bs v^{\ast})}$ (see, \cite{Mum76}, Theorem $2.33$).
In particular, $G$-polystability implies $G$-semistability, since a closed $G$-orbit cannot contain the origin.
\end{remark}
%From the above definition, one can see that $G$-polystability implies $G$-semistability because $G$-orbit itself never contains the origin. 

\vskip 7pt

We now recall the definition of the Chow form of an irreducible complex projective variety; see \cite{GKZ94} for details. Let $X\hookrightarrow \PP^N$ be an $n$-dimensional irreducible
complex projective variety of degree $d\geq 2$. Recall that the
Grassmannian $\mathbb{G}(k,\PP^N)$ parameterizes $k$-dimensional
projective linear subspaces of $\PP^N$.
\begin{definition}\rm
The \emph{associated hypersurface} of $X\rightarrow \PP^N$ is the subvariety 
%in $\mathbb{G}(N-n-1, \PP^N)$ which is given by
\[
Z_X:=\set{L\in \GG(N-n-1,\PP^N)| L\cap X\neq \varnothing}.
\]
\end{definition}
The fundamental properties of $Z_X$ are as follows (see \cite[p. $99$]{GKZ94}):
\begin{enumerate}
\item $Z_X$ is irreducible;
\item $\mathrm{codim\:} Z_X=1$, i.e., $Z_X$ is a divisor in
$\GG(N-n-1,\PP^N)$;
\item $\deg Z_X=d$ with respect to the Pl\"ucker embedding; 
\item $Z_X$ is defined by the vanishing of a section 
\[
R_X^{\ast}\in H^0\bigl(\GG(N-n-1,\PP^N),\mathcal O(d)\bigr).
\]
\end{enumerate}
The section $R_X^{\ast}$ is called the \emph{Chow form} of $X$, and is uniquely %. Note that $R_X^{\ast}$ can be 
determined up to a non-zero scalar multiple. %Setting
Let
\[
V:=H^0\bigl(\GG(N-n-1,\PP^N),\mathcal O(d)\bigr), \qquad R_X \in \PP(V)
\]
be the corresponding point.
%which is the projectivization of $R_X^{\ast}$, 
We call $R_X$ the \emph{Chow point} of $X$. Since $G={{SL}}(N+1,\CC)$ acts naturally on $\PP(V)$, 
we can define GIT-stability of $R_X$ with respect to this action.
%we can define ${\rm{SL}}(N+1)$-polystability (resp. semistability) of $R_X$.
%Def
\begin{definition}\rm
Let $X\subset \PP^N$ be an irreducible $n$-dimensional complex projective variety. 
We say that $X$ is \emph{Chow polystable (resp. Chow semistable)} if its Chow point $R_X$ %of $X$ 
is ${{SL}}(N+1,\CC)$-polystable (resp. semistable).
If $X$ is not Chow semistable, then $X$ is said to be {\emph{Chow unstable}}.
\end{definition}
%Def
\begin{definition}\rm
Let $(X,L)$ be a polarized variety, and let 
\[
\Psi_k: X \hookrightarrow \PP\left(H^0(X;\mathcal O_X(L^k))\right)
\] 
be the Kodaira embedding.
We say that $(X,L)$ is {\emph{asymptotically Chow polystable}} (resp. {\emph{asymptotically Chow semistable}}) if 
$\Psi_k(X)\subset \PP(H^0(X;\mathcal O_X(L^k)))$ is Chow polystable (resp. Chow semistable) for all sufficiently large integers $k\gg 0$.
\end{definition}

Next, we recall some basic facts from toric geometry; see \cite{CLS11} for details.  
A {\emph{toric variety}} $X$ is a normal algebraic variety equipped with an effective holomorphic action of the algebraic torus $T_\CC:=(\CC^\times)^n$, such that $\dim_\CC X=n$.
Let $T_\bbR:=(S^1)^n$ be the compact real torus, %in $T_\CC$ 
and $\mathfrak{t}_\bbR$ denote its Lie algebra. Set 
\[
N_\bbR:=J\mathfrak{t}_\bbR\cong \bbR^n,
\] %be the associated Euclidean space 
where $J$ denotes the complex structure on $T_\CC$. 
%We denote the dual space $\Hom (N_\bbR, \bbR)\cong \bbR^n$ of $N_\bbR$ by $M_\bbR$.
Let 
\[
M_\bbR:=\Hom (N_\bbR, \bbR)\cong \bbR^n
\]
be the dual vector space. If $M$ denotes the lattice of algebraic characters of $T_\CC$, then
%Setting the group of algebraic characters of $T_\CC$ by $M$, we see that 
$M_\bbR=M\otimes_\bbZ\bbR$, and the dual lattice $N:=\Hom(M,\bbZ)$ consists of one-parameter subgroups of $T_\CC$.
%Let us denote the dual lattice of $M$ by $N$. Then, $N$ consists of the algebraic one parameter subgroups of $T_\CC$.

Let $\Delta \subset M_\bbR\cong \bbR^{n}$ be an $n$-dimensional integral convex polytope. 
Define the cone over $\Delta\times\set{1}$ by
 \[
 C(\Delta):=\set{(\bs x, r)\in\mathbb{R}^n\times\mathbb{R}\mid r>0,\, r^{-1}\bs x\in \Delta}\cup\set{\bs 0}.
 \]
Note that $C(\Delta)$ is a strongly convex cone with vertex at the origin.
%to be the cone over $\Delta\times\set{1}$ with the vertex at the origin. 
By Gordan's Lemma, the semigroup $S_\Delta:= C(\Delta)\cap\mathbb{Z}^{n+1}$ is finitely generated.
%by Gordan's lemma.
Let $\mathbb{C}[S_{\Delta}]$ denote its semigroup algebra. The character corresponding to $(\bs m,k)\in S_\Delta$ is 
denoted by $\chi^{\bs m}t^k$, and $\mathbb{C}[S_\Delta]$ is naturally graded by the height: $\deg(\chi^{\bs m}t^k)=k$. Thus we obtain a graded $\mathbb{C}$-algebra
\[
\mathbb{C}[S_\Delta]=\smalloplus_{k=0}^{+\infty}R_k, \qquad R_k:=\Set{f\in \mathbb{C}[S_\Delta]\mid\deg f=k}.
\]
%from the polytope $\Delta$. 
We define the polarized toric variety associated to $\Delta$ by %$(X_\Delta, L_\Delta)$ by
\[
(X_\Delta, L_\Delta):=\bigl(\operatorname{Proj}\,(\mathbb{C}[S_\Delta]), \mathcal{O}_{X_\Delta}(1)\bigr).
\]
It is well-known that $X_\Delta$ is smooth ad projective if and only if $\Delta$ is Delzant.

%Sec2.2
\subsection{Simplex triangulations}\label{sec:SimpTri}
%Def
\begin{definition}\label{def:isom}\rm
Let $\Delta, \Delta' \subset \bbR^n$ be integral polytopes. We say that $\Delta$ is \emph{isomorphic} to $\Delta'$ if there exist translations $T_1$, $T_2$ and an element $g \in SL(n,\bbZ)$ such that, for every $k\in \bbN$,
\[
k\Delta'=T_2\circ g \circ T_1 (k\Delta),
\] 
and 
\[
\bs p\in k\Delta \cap \bbZ^n \qquad \text{if and only if} \qquad T_2\circ g \circ T_1(\bs p) \in k\Delta' \cap \bbZ^n.
\]  
\end{definition}
Note that $\Delta$ and $\Delta'$ are not assumed to be full-dimensional.
%Def
\begin{definition}\rm
A $k$-dimensional polytope $S$ is called a {\emph{$k$-dimensional enlarged simplex}} if it is isomorphic to 
$\conv\Set{\bs e_0,r \bs e_1,r\bs e_2, \dots, r\bs e_k}$ for some $r\in \bbN$, where
\[
\bs e_0=(0, \dots ,0), \quad \bs e_1=(1,0, \dots, 0), \quad \dots ,\quad \bs e_k=(0,\dots, 0, 1).
\]
We call $S$ a {\emph{unimodular simplex}} if it is isomorphic to 
$\conv\Set{\bs e_0, \bs e_1, \dots, \bs e_k}$.
\end{definition}
We now introduce several special classes of triangulations of integral polytopes.
%Def
\begin{definition}\label{def:SimpTri}\rm
Let $\Delta$ be an integral polytope.
\begin{enumerate}
\item An {\emph{integral triangulation}} of $\Delta$ is a triangulation in which all vertices of every simplex belongs to $\bbZ^n$.   
\item A {\emph{simplex triangulation}} of $\Delta$ is an integral triangulation in which every simplex is unimodular (i.e., isomorphic to the standard simplex).
 \end{enumerate} 
\end{definition}
%Rem
\begin{remark}\rm
From this point onward, we assume that every simplex is {\emph{unimodular}}, unless stated otherwise. 
In particular, when we use the term triangle, we allow the possibility that it is {\emph{non-unimodular}}.
\end{remark}

In particular, for any simplex triangulation of $\Delta$ and any lattice point $\bs p\in \Delta \cap \mathbb Z^n$,
there exists at least one simplex $S$ having $\bs p$ as a vertex.

Let $H=\set{\bs x\in \bbR^n | h(\bs x)=0}$ be a hyperplane in $\bbR^n$.
We define the associated upper half-space by
\[
H^+:=\set{\bs x \in \bbR^n | h(\bs x)\geq 0}.
\]
Let $\bs v_1, \dots ,\bs v_r$ be vectors contained in a half-space $H^+$ for some hyperplane $H$.
For a point $\bs p \in \bbR^n$, we define the {\emph{infinite half-cone}} with vertex $\bs p$ by
\[
C(\bs p; \bs v_1, \dots, \bs v_r):=\Set{\bs q \in \bbR^n | \bs q=\bs p+t_1\bs v_1 +\dots + t_r\bs v_r, ~ t_i \geq 0}.
\]
%We say that $C(p;v_1, \dots,v_r)$ is {\emph{an infinite half cone of the vertex $p$ with generators $\set{v_1, \dots ,v_r}$}} if the following holds: 
%for each $q\in C(p;v_1, \dots,v_r)$, there exist $t_1, \dots ,t_r\geq 0$ such that 
%\[
%q=p+t_1v_1+\cdots +t_rv_r.
%\] 
%We will denote the cone $C(p;v_1, \dots,v_r)$ as $C(p)$ if the generators are known.
When the generators are clear from the context, we simply write $C(\bs p)$.

%For any integral polytope $\Delta$, let $p_1, \dots ,p_R$ be the vertices of $\Delta$. 
Let $\Delta$ be an integral polytope with vertices $\bs p_1, \dots , \bs p_R$.
For each $i$, let $C(\bs p_i)$ be the strongly convex polyhedral cone generated by the edges emanating from 
$\bs p_i$. Then, $\Delta$ can be written as the intersection of these cones:
\begin{equation}\label{eq:decomp}
\Delta=\smallcap_{i=1}^RC(\bs p_i)
\end{equation}
%in the form of the intersection of $C(p_i)$.
Moreover, for a sufficiently small neighborhood $U_{\bs p_i}$ of $\bs p_i$, we have %the equality
\[
C(\bs p_i)\cap U_{\bs p_i}=\Delta \cap U_{p_i}. 
 \] 
%Note that if 
%\[
%\Delta=\bigcap_{i=1}^RC(p_i)
%\] 
%holds, then we have
For any positive integer $k$, $\eqref{eq:decomp}$ implies that
\[
k\Delta=\smallcap_{i=1}^R C(k \bs p_i),
\] 
where 
\[
C(\bs p_i)=C(k \bs p_i)-(k-1)\bs p_i:=\set{\bs x-(k-1)\bs p_i | \bs x\in C(k\bs p_i)}.
\] 
%Hence, each cone shares the same triangulation.
In particular, each cone $C(k \bs p_i)$ admits a compatible triangulation.
We denote a simplex triangulation of $C(k\bs p_i)$ by $\mathcal T(C(k \bs p_i))$.

Let $\ell:\bbR^n \rightarrow \bbR$ be a non-zero linear function and let $c\in \bbR$.
%which is not identically zero, and a real constant $c\in \bbR$, the set
The set
\[
H^-=\Set{\bs x\in \bbR^n | \ell(\bs x)\leq c}
\]
is called a (lower) \emph{half-space}. A \emph{polyhedron} is an intersection of finitely many 
half-spaces. It is called \emph{integral} if all of its vertices belong to $\bbZ^n$.
 %Def
\begin{definition}\label{def:n_ik}\rm
 Let $\mathcal{B}_k:=\set{B\subset k\Delta | B\text{ is an integral  polyhedron}}$.
Following \cite{Lee25}, we define a function
\[
n_{i,k}: (C(k \bs p_i)\cap \bbZ^n)\times \mathcal{B}_k \longrightarrow \bbZ
\]
by
\[
n_{i,k}(\bs q,B):= \#\set{S|\text{$S$ is a simplex in $B$ containing $\bs q$}}.
\]
%$n_{i,k}$ is a function defined by
%\begin{equation*}
%\xymatrix@R=-1ex{
%n_{i,k}: (C(kp_i)\cap \bbZ^n)\times \mathcal{B}_k \ar[r]& \bbZ\\
%\quad \! \din &\din\\
 %(q,B) \ar@{|->}[r]& \#\set{S|\text{$S$ is a simplex in $B$ touching $q$}}.
 %}
%\end{equation*}
%If we take $B$ as $k\Delta$, then $n_{i,k}(q,B)$ is denoted by $n_{i,k}(q)$: 
In particular, when $B=k\Delta$, we write
\begin{equation}\label{eq:n_ik}
n_{i,k}(\bs q):= n_{i,k}(\bs q,k\Delta).
\end{equation}
\end{definition}
%Obviously, the following equality holds by Definition \ref{def:n_ik}.
The following property is immediate from the definition.
%Lem
\begin{lemma}\label{lemm partition estimate}
Let $B$ be partitioned into polyhedra $B_1, \dots ,B_{N}$. Then, for any vertex $\bs p_i\in \mathcal V(\Delta)$ and any
lattice point $\bs q \in C(k \bs p_i)\cap \bbZ^n$, we have
\[
\sum_{r=1}^Nn_{i,k}(\bs q,B_{r})=n_{i,k}(\bs q,B).
\]
%As a consequence, we see that
In particular, 
\[
n_{i,k}(\bs q,B)\leq n_{i,k}(\bs q).
\]
\end{lemma}

%Sec2.3
\subsection{Partition of unity}\label{sec:partition}
The goal of this subsection is to compare the integral $\int_{k\Delta}f(\bs x){d\mu}$ with the lattice-point summand 
$\sum_{\bs p\in k\Delta\cap \bbZ^n}f(\bs p)$, using only local data associated with the vertex cones $C(\bs p_i)$.
To this end, we construct a partition of unity $\set{\delta_i}$ on $\Delta$.
%Sec2.3.1
\subsubsection{Construction of the partition of unity}
Let $\Delta$ be an integral polytope %associated with a polarized toric variety, 
and let $C(\bs p_i)$ denote the cone at a vertex $\bs p_i$ as described in Section \ref{sec:SimpTri}.
%Let $\delta_i:\Delta \rightarrow \bbR$ be a partition of unity of $\Delta$ constructed by the following:
We construct a partition of unity $\delta_i:\Delta \rightarrow \bbR$ as follows.

\begin{enumerate}
\item[{\emph{Step 1}}.] For any $\bs p \in \Delta$ and $t>0$, define %we can define $t\Delta(p)$ by 
\[
t\Delta(\bs p):=\Set{s(\bs x- \bs p)+ \bs p | 0\leq s\leq t, ~ \bs x\in \Delta }.
\] 
For each vertex $\bs p_i$ and each $t>0$, choose a smooth function $\eta_{i,t}$ satisfying
%we have a smooth function $\eta_{i,t}$ such that 
\begin{itemize}
\item $\supp(\eta_{i,t})\subset \left(1-\frac{t}{2}\right)\Delta(\bs p_i)$, and
\item $\eta_{i,t}(\bs x)>0$ for all $\bs x\in (1-t)\Delta(\bs p_i)$.
\end{itemize}
\item[{\emph{Step 2}}.] Choose $t>0$ sufficiently small so that 
\[
\Delta=\smallcup_{i=1}^R(1-t)\Delta(\bs p_i).
\] 
Define 
\[
\eta(\bs x):=\sum_{i=1}^R\eta_{i,t}(\bs x).
\] 
By construction, $\eta(\bs x)>0$ for every $\bs x\in \Delta$.
\item[{\emph{Step 3}}.] Define $\delta_i(\bs x):={\eta_{i,t}(\bs x)}/{\eta(\bs x)}$. Then each $\delta_i(\bs x)$ is smooth and 
\[
\sum_{i=1}^R\delta_i(\bs x)=1, %\qquad \text{for all} 
\qquad \bs x\in \Delta.
\]
\end{enumerate}
Moreover, for sufficiently small neighborhoods $B_{\vep}(\bs p_j)$ of the vertices $\bs p_j$, we have 
\[
\delta_i(\bs x)=\delta_{ij}, \qquad %\text{for}\quad  
\bs x\in  B_{\vep}(\bs p_j),
\]  
where $\delta_{ij}$ denotes the Kronecker delta. 
We extend $\delta_i$ to $k \Delta$ by defining
\[
\delta_{i,k}(\bs x):=\delta_i\left(\dfrac{\bs x}{k}\right), \qquad \bs x \in k\Delta,
\]
and setting $\delta_{i,k}(\bs x)=0$ for $\bs x \in C(k \bs p_i)\backslash k\Delta$.

Since $\Delta$ is compact, the functions $\delta_i$ can be chosen so that there exists a constant
$c>0$ satisfying
\[
\norm{D^{\bs \alpha}\delta_i(\bs x)}<c 
\]
for all multi-indices $\bs \alpha$ with $\norm{\bs \alpha}\le 2$.
Here, for a smooth function $f(\bs x)$ on $\bbR^n$ and a multi-index $\bs \alpha=(\alpha_1,\alpha_2, \dots, \alpha_n)$ with $|\bs \alpha|=\alpha_1+\alpha_2+\dots +\alpha_n$, we write
\[
D^{\bs \alpha}f:=\frac{\partial^{|\bs \alpha|}f}{\dd x_{1}^{\alpha_1}\dots \dd x_{n}^{\alpha_n}}.
\]
Next, we estimate integrals of the form
\[
\int_{S}\delta_{i,k}(\bs x)f(\bs x) d\mu,
\] 
where $S$ is a simplex. To this end, we first analyze integrals over the standard simplex
\[
\Delta_n:=\conv\Set{\bs e_0, \bs e_1, \dots, \bs e_n} \subset \mathbb R^n,
\] 
where
\[
\bs e_0=(0, \dots, 0), \quad  \bs e_1=(1,0, \dots, 0), \quad \dots , \quad  \bs e_n=(0,\dots, 0,1).
\]
To proceed, we recall the Bramble-Hilbert lemma (see \cite{Man07} for details).

Let $\Omega\subset \bbR^n$ be a convex domain, and let $W_p^m(\Omega)$ denote the Sobolev space consisting of 
functions $u$ on $\Omega$ whose weak derivatives $D^{\bs \alpha} u$ of order $|\bs \alpha| \leq m$ belong to $L^p(\Omega)$. The Sobolev seminorm on $W_p^m(\Omega)$ is defined by
\[
|u|_{W_{p}^m(\Omega)}:=\left(\sum_{|\bs \alpha|=m}\|D^{\bs \alpha}u\|_{L^p(\Omega)}^p\right)^{\! {1}/{p}}, \qquad 
1\leq p <\infty,
\]
and
\[
|u|_{W_{\infty}^m(\Omega)}:=\max_{|\bs \alpha|=m}\|D^{\bs \alpha}u\|_{L^{\infty}(\Omega)}. 
\]

%Lem:BrHil
\begin{lemma}[Bramble-Hilbert Lemma (\cite{BrHi70}, Theorem 2)]\label{lem:BH} 
Let $\ell$ be a bounded continuous linear functional on $W_p^m(\Omega)$, and let $\|\ell\|_{W_p^m(\Omega)'}$ denote its dual norm.
Assume that $\ell(g)=0$ for every polynomial $g$ of degree at most $m-1$.
Then there exists a constant $C=C(\Omega)$ such that, for all $u\in W_p^m(\Omega)$,
\[
|\ell(u)|\leq C \|\ell\|_{W_p^m(\Omega)'} |u|_{W_p^m(\Omega)}.
\] 
\end{lemma} 

%Lem
\begin{lemma}\label{lem: deeper estimate}
Let $S$ be an $n$-dimensional unimodular simplex with vertices $\bs p_1,\dots, \bs p_{n+1}$. 
Set $\bs v_{ij}:=\bs p_j-\bs p_i$ and $f_{\bs v_{ij}}:=\langle\nabla f, \bs v_{ij}\rangle$.
Then, for any affine functions $a(\bs x), f(\bs x):S\to \mathbb R$, we have
\[
\int_S a(\bs x)f(\bs x)d\mu=\dfrac{1}{(n+1)!}\sum_{i=1}^{n+1}a(\bs p_i)f(\bs p_i)+\dfrac{1}{(n+2)!}\sum_{i=1}^{n+1}\sum_{j=1\atop j\ne i}^{n+1}a(\bs p_i)f_{\bs v_{ij}}(\bs p_i).
\] 
\end{lemma}
\begin{proof}
By an affine unimodular change of coordinates, 
we may assume that 
\[
S=\conv\Set{ \bs e_0,\dots, \bs e_n}.
\] 
Both sides of the identity are bilinear in $a$ and $f$. Since the space of affine functions on $S$ is generated by
$1,x_1, \dots, x_n$, it suffices to verify the formula for the following cases:
\begin{enumerate}
     \item $a(\bs x)=f(\bs x)=1$;
          \item $a(\bs x)=x_i$, $f(\bs x)=1$;
     \item $a(\bs x)=1$, $f(\bs x)=x_i$;
        \item $a(\bs x)=f(\bs x)=x_i$;         and
        \item $a(\bs x)=x_i$, $f(\bs x)=x_j$ for $i\neq j$;
\end{enumerate} 
By symmetry, it suffices to consider representative cases $x_1$ and $x_2$.
\vskip 5pt
\noindent (1) Case $a(\bs x)=f(\bs x)=1$. We have
\[
\int_S a(\bs x)f(\bs x)d\mu=\vol(S)=\dfrac{1}{n!}.
\] 
On the other hand, 
\begin{align*}
\dfrac{1}{(n+1)!}\sum_{i=1}^{n+1}a(\bs p_i)f(\bs p_i)+\dfrac{1}{(n+2)!}\sum_{i=1}^{n+1}\sum_{j=1\atop j\neq i}^{n+1}a(\bs p_i)f_{\bs v_{ij}}(\bs p_i) 
= \dfrac{n+1}{(n+1)!}
=\dfrac{1}{n!},
\end{align*}
as $f_{\bs v_{ij}}(\bs p_i)=0$. Hence the identity holds.
\vskip 5pt
%%%%%%%%%%%%%%%%%%
\noindent (2) Case $a(\bs x)=x_1$, $f(\bs x)=1$. We compute
\[
\int_S a(\bs x)f(\bs x)d\mu=\int_S x_1d\mu=\vol(\Delta_{n+1})=\dfrac{1}{(n+1)!},
\] 
where $\Delta_{n+1}$ is the $(n+1)$-dimensional standard simplex.
Since $f$ is constant, $\nabla f=0$, and hence all directional derivatives vanish. Moreover,
$a(\bs p_i)=0$ unless $\bs p_i=\bs e_1 $, where $a(\bs e_1)=1$. Thus,
\begin{align*}
\dfrac{1}{(n+1)!}\sum_{i=1}^{n+1}a(\bs p_i)f(\bs p_i)=\frac{1}{(n+1)!},
    \end{align*}
    and the second term is zero. The identity follows.
\vskip 5pt
%%%%%%%%%%%%%%%%%%
\noindent (3)
Case $a(\bs x)=1$, $f(\bs x)=x_1$. Similarly,  
\[
\int_S a(\bs x)f(\bs x)d\mu=\int_S x_1d\mu=\vol(\Delta_{n+1})=\dfrac{1}{(n+1)!}.
\] 
Now $a(\bs p_i)=1$ for all $i$, and 
\[
\nabla f=(1,0,\dots, 0).
\]
Thus, 
\begin{equation}\label{eq:DirDeri}
\braket{\nabla f, \bs v_{ij}}=
\begin{cases}
1 & \text{if } \quad j=1, i\neq 1,\\
-1 & \text{if }\quad  i=1, j\neq 1,\\
0 & \text{otherwise}.
\end{cases}
\end{equation}
Summing over $i,j$, one obtains that
\begin{align*}
&\dfrac{1}{(n+1)!}\sum_{i=1}^{n+1}a(\bs p_i)f(\bs p_i)+\dfrac{1}{(n+2)!}\sum_{i=1}^{n+1}\sum_{j=1 \atop j\neq i}^{n+1}a(\bs p_i)f_{\bs v_{ij}}(\bs p_i)\\
=& \dfrac{1}{(n+1)!}+\dfrac{1}{(n+2)!}\left\{\sum_{j=2}^n f_{\bs v_{1j}}(\bs p_i)+\sum_{i=2}^n f_{\bs v_{i1}}(\bs p_i)
\right\} 
=\dfrac{1}{(n+1)!}.
    \end{align*}
Hence the identity holds.
\vskip 5pt
%%%%%%%%%%%%%%%%%%
\noindent (4) Case $a(\bs x)=f(\bs x)=x_1$. A direct computation on the standard simplex
\[
S=\conv\Set{\bs e_0, \bs e_1, \dots, \bs e_n}=\Delta_{n}
\]
gives
\begin{align}
%\begin{split}\label{eq:Int_case4}
\int_S x_1^2d\mu =&\int_0^1x_1^2 dx_1 \int_0^{1-x_1}dx_2 \cdots \int_0^{1-x_1-\cdots-x_{n-1}}dx_n \notag \\
=&\int_0^1x_1^2 dx_1 \int_0^{1-x_1} dx_2 \cdots \int_0^{1-x_1-\cdots-x_{n-2}}(1-x_1-\cdots-x_{n-1})dx_{n-1} \label{eq:Int_case4} \\
=& \cdots =\int_0^1\dfrac{x_1^2(1-x_1)^{n-1}}{(n-1)!} dx_1=\frac{2}{(n+2)!}. \notag
\end{align}
On the other hand, only the vertex $\bs e_1$ contributes to the first term:
\[
\frac{1}{(n+1)!}a(\bs e_1)f(\bs e_1)=\frac{1}{(n+1)!}.
\]
For the second term, using the same computation of directional derivatives in $\eqref{eq:DirDeri}$, 
\begin{align*}
\sum_{i, j\neq i}a(\bs p_i)\braket{\nabla f(\bs p_i), \bs v_{ij}}=\sum_{j=2}^{n+1}a(\bs e_1)\cdot(-1)=-n,
\end{align*}
where we used $a(\bs p_i)f_{\bs v_{ij}}(\bs p_i)=0$ unless $\bs p_i =\bs e_1$ in the above computation.
Thus,
\begin{align*}
&\frac{1}{(n+1)!}\sum_{i=1}^{n+1}a(\bs p_i)f(\bs p_i)+\frac{1}{(n+2)!}\sum_{i=1}^{n+1}\sum_{j=1 \atop j \neq i}^{n+1}a(\bs p_i)f_{\bs v_{ij}}(\bs p_i)\\
=&\frac{1}{(n+1)!}-\frac{n}{(n+2)!}=\frac{2}{(n+2)!},
\end{align*} 
which agrees with the integral $\eqref{eq:Int_case4}$.

\vskip 5pt
%%%%%%%%%%%%%%%%%%
\noindent (5)
Case $a(\bs x)=x_1$, $f(\bs x)=x_2$. A direct computation yields
\begin{align*}
\int_S x_1x_2d\mu=&\int_0^1x_1\, dx_1\int_0^{1-x_1}x_2\, dx_2 \int_0^{1-x_1-x_2}dx_3\,  \dots \int_0^{1-\sum_{i=1}^{n-1}x_i}dx_n \\
=&\int_0^1x_1 \left( \int_0^{1-x_1}
      \dfrac{x_2(1-x_1-x_2)^{n-2}}{(n-2)!}dx_2 \right)dx_1\\
=&\frac{1}{(n-2)!}\int_0^1 x_1 \left( \int_0^{1-x_1}(1-x_1-u)u^{n-2}du\right) dx_1\\
      =&\frac{1}{(n-2)!} \int_0^1x_1\left[\dfrac{(1-x_1)(1-x_1)^{n-1}}{n-1}-\dfrac{(1-x_1)^n}{n}\right]dx_1\\
      =&\frac{1}{n!}\int_0^1x_1(1-x_1)^{n}dx_1=\dfrac{1}{(n+2)!}.
  \end{align*}
Since no vertex of $S=\Delta_n$ simultaneously satisfies $x_1=x_2=1$, the first term vanishes;
\[
\dfrac{1}{(n+1)!}\sum_{i=1}^{n+1}a(\bs p_i)f(\bs p_i)=0.
\]
For the second term, only $\bs p_i=\bs e_1$, $j=2$ contributes, giving
\begin{align*}
\dfrac{1}{(n+2)!}\sum_{i=1}^{n+1}\sum_{j=1 \atop j\ne i}^{n+1}a(\bs p_i)f_{\bs v_{ij}}(\bs p_i)
=&\dfrac{1}{(n+2)!} a(\bs e_1)f_{\bs v_{12}}(\bs e_1)
=\dfrac{1}{(n+2)!}.
\end{align*}
Hence the identity holds.

\vskip 5pt

Since both sides are bilinear in $a$ and $f$, and the above cases form a basis for affine functions, the identity follows.
\end{proof}
We now apply the preceding results to define a family of linear functionals.
Let $S$ be an $n$-dimensional unimodular simplex with vertices are $\bs p_1,\dots, \bs p_{n+1}$,
and let $a(\bs x)$ be an affine function on $S$. For a smooth function $f$, define
\begin{equation}\label{def:LinFunL_a}
L_a(f):=\int_S a(\bs x)f(\bs x)d\mu-\dfrac{1}{(n+1)!}\sum_{i=1}^{n+1}a(\bs p_i)f(\bs p_i)-\dfrac{1}{(n+2)!}\sum_{i=1}^{n+1}\sum_{j=1 \atop j \neq i}^{n+1}a(\bs p_i)f_{\bs v_{ij}}(\bs p_i),
\end{equation}
where 
\[
f_{\bs v_{ij}}:=\langle\nabla f, \bs v_{ij}\rangle, \qquad \bs v_{ij}=\bs p_j -\bs p_i
\]
Lemma \ref{lem: deeper estimate} together with the Bramble-Hilbert Lemma immediately yields the following estimate.
%Cor
\begin{corollary}\label{coro: improved local estimate}
For any affine function $a(\bs x)$ and any smooth function $f(\bs x)$, there exists a constant $C_S>0$ such that 
\begin{align*}
&\left|\int_S a(\bs x)f(\bs x)d\mu-\frac{1}{(n+1)!}\sum_{i=1}^{n+1}a(\bs p_i)f(\bs p_i)
- \frac{1}{(n+2)!}\sum_{i=1}^{n+1}\sum_{j\neq i}a(\bs p_i)f_{\bs v_{ij}}(\bs p_i)\right|\\
   \leq& C_S\left(\int_S|a(\bs x)|d\mu\right)  \sup_{\bs x\in S \atop |\bs \alpha|=2}\left|D^{\bs \alpha} f(\bs x)\right|.
\end{align*}
\end{corollary}
\begin{proof}
For a fixed affine function $a(\bs x)$, let $L_a(f)$ be the linear functional defined in $\eqref{def:LinFunL_a}$.
Then $L_a$ is a bounded linear functional on $W_\infty^2(S)$, since $S$ is bounded. 
By Lemma \ref{lem: deeper estimate}, $L_a(f)=0$ for all polynomials of degree at most $1$. Applying the Bramble-Hilbert lemma (Lemma \ref{lem:BH}), we obtain  
\[
L_a(f)\leq C_S\norm{L_a}_{W_\infty^2(S)'} \sup_{\bs x\in S \atop \norm{\bs \alpha}=2}\left|D^{\bs \alpha} f(\bs x)\right|.
\]
Since $\norm{L_a}_{W_\infty^2(S)'}$ is bounded by a constant multiple of $\int_S|a(\bs x)|d\mu$, 
the desired estimate follows.
\end{proof}
%Rem
\begin{remark}\rm
The constant $C_S$ depends on the geometry of the simplex $S$, in particular on its diameter and shape.
Therefore, when applying Corollary \ref{coro: improved local estimate} to triangulations of $k\Delta$,
it is important that only finitely many simplex types occur. Under this assumption, the constants $C_S$ may be chosen uniformly over all simplices in the triangulation.
\end{remark}

%Sec2.4
\subsection{Type F triangulation}
%Def
\begin{definition}\label{def:TypeF}\rm
Let $C(\bs p_\alpha)$ be the vertex cone at $\bs p_\alpha$.
A \emph{type F triangulation} of $C(\bs p_\alpha)$ is a simplex triangulation $\mathcal T(C(\bs p_\alpha))$ satisfying the following conditions:
\begin{enumerate}
\item {\bf{Scaling invariance:}} For any positive integer $k$, 
\[
\mathcal T(C(\bs p_\alpha))=\mathcal T(C(k\bs p_\alpha)).
\]
\item {\bf{Finiteness of shapes:}} There exist finitely many simplices $S_1, \dots , S_M$ such that for every simplex 
$S\in \mathcal T(C(\bs p_\alpha))$, there exist a lattice vector $\bs c\in \bbZ^n$ and an index $j\in \set{1, \dots, M}$ such that 
\[
S-\bs c=S_j.
\] 
\end{enumerate}
In other words, every simplex in the triangulation is a lattice translate of one of finitely many {\emph{model simplices}}.
\end{definition}

%Cor
\begin{corollary}\label{coro: global estimate refine}
Let $\Delta=\conv\set{\bs p_{\alpha}}_{\alpha=1}^R$ be an integral polytope with vertices $\bs p_{\alpha}$.
Let $f$ be a non-negative convex function on $k\Delta$, and let $\set{\delta_{\alpha}}_{\alpha=1}^R$ be a partition of unity as constructed in Section \ref{sec:partition}. Assume that, for each vertex $\bs p_\alpha$, the cone $C(k\bs p_{\alpha})$ admits a type F triangulation
$\mathcal T:=\mathcal T(C(k\bs p_\alpha))$. 
Then 
\begin{align*}
(n+1)!\int_{k\Delta}  f(\bs x){d\mu} &\leq  \sum_{\bs p\in k\Delta\cap \bbZ^n}n_{k}(\bs p)f(\bs p)
 %\frac{n_{k}(\bs p)f(\bs p)}{(n+1)!} 
+\dfrac{C}{k^2}\sum_{\alpha=1}^R
\sum_{\bs p\in k\Delta\cap \bbZ^n}n_{\alpha,k}(\bs p)f(\bs p)\\
%\dfrac{n_{\alpha,k}(\bs p)f(\bs p)}{(n+1)!}\\
& +\frac{1}{(n+2)}\sum_{\alpha=1}^R\sum_{\bs p\in k\Delta\cap \bbZ^n}\sum_{S\in \mathcal T  \atop \bs p\in \mathcal{V}(S)}
\sum_{i=1}^{n+1}\sum_{j=1\atop j\neq i}^{n+1}f(\bs p)\delta_{\alpha,k, \bs v_{ij}^S}(\bs p),
\end{align*}
where:
\begin{itemize}
\item $\delta_{\alpha,k}(\bs x):=\delta_{\alpha}(\bs x/k)$,  
\item $\delta_{\alpha,k, \bs v_{ij}^S}:=
\braket{\nabla \delta_{\alpha,k},\bs v_{ij}^S}$, with $\bs v_{ij}^S:=\bs p_j - \bs p_i$ for $\bs p_i, \bs p_j \in \mathcal V(S)$,
\item $n_k(\bs p):=\max\set{n_{\alpha,k}(\bs p) | \alpha=1,\dots, R}$,
\item and the constant $C$ is given by 
\begin{equation}\label{def:constC}
C:=\max_{1\leq \alpha \leq R}C_\alpha, \quad \text{where} \quad
C_\alpha:=\max_{S\in \mathcal T(C(k\bs p_{\alpha}))}C_S\cdot \max_{\bs x\in \Delta}\|\mathrm{Hess}(\delta_{\alpha}(\bs x))\|.
\end{equation}
\end{itemize}
Here the sum over $S\in \mathcal T$,  $\bs p\in \mathcal{V}(S)$ runs over all simplices $S\in \mathcal T$ having the lattice point $\bs p \in k\Delta \cap \bbZ^n$ as a vertex.
\end{corollary}
%proof
\begin{proof}
Let $f:k\Delta\rightarrow \bbR$ be a convex function. Fix $\alpha$, and let $\mathcal T:=\mathcal T(C(k \bs p_{\alpha}))$ be a type F triangulation of the cone $C(k \bs p_{\alpha})$.
Let $a_f$ denote the piecewise-linear interpolation of $f$ with respect to the triangulation $\mathcal T(C(k \bs p_{\alpha}))$; namely, $a_f$ is affine on each simplex $S\in \mathcal T(C(k \bs p_{\alpha}))$ and satisfies
\[
a_f(\bs p)=f(\bs p)
\] 
for every lattice point $\bs p\in k\Delta\cap \bbZ^n$.
Applying Corollary \ref{coro: improved local estimate} to each simplex $S\in \mathcal T$, with $a(\bs x)=a_f(\bs x)$ and $f(\bs x)=\delta_{\alpha, k}(\bs x)$, we obtain
\begin{align*}
 &\int_{k\Delta}  f(\bs x)\delta_{\alpha}\left(\dfrac{\bs x}{k}\right){d\mu}
 \leq \int_{k\Delta}  a_f(\bs x)\delta_{\alpha}\left(\dfrac{\bs x}{k}\right){d\mu}\\
 \leq&  \frac{1}{(n+1)!}\sum_{\bs p\in k\Delta \cap \bbZ^n}n_{\alpha,k}(\bs p)a_f(\bs p)\delta_{\alpha}\left(\dfrac{\bs x}{k}\right)
 +\dfrac{C}{k^2}\int_{k\Delta}a_f(\bs x)\delta_{\alpha}\left(\dfrac{\bs x}{k}\right){d\mu} \\
& + \frac{1}{(n+2)!}\sum_{\bs p\in k\Delta\cap \bbZ^n}\sum_{S\in \mathcal T \atop \bs p\in \mathcal{V}(S)}\sum_{i, j\ne i}a_f(\bs p)\delta_{\alpha,k, \bs v_{ij}^S}(\bs p),
\end{align*}
for the constant $C$ defined in $\eqref{def:constC}$.
%\[
%C=\max_{S\in T}C_S\cdot \max_{\bs x\in \Delta}\|\mathrm{Hess}(\delta_{\alpha}(\bs x))\|,
%\]
Here we used the scaling estimate
\[
\max_{|\bs \beta|=2}\left|D^{\bs \beta}\delta_{ \alpha}\left(\dfrac{\bs x}{k}\right)\right|=\dfrac{1}{k^2}\max_{|\bs \beta|=2}|D^{\bs \beta}\delta_{\alpha}(\bs x)|=\dfrac{1}{k^2}\|\mathrm{Hess}(\delta_{\alpha}(\bs x))\|
\]
in the above inequality. By the definition of $a_f$, we have $a_f(\bs p)=f(\bs p)$ for lattice points $\bs p$, hence
we see that
\begin{align*}
(n+1)!&\int_{k\Delta}  f(\bs x)\delta_{\alpha}\left(\dfrac{\bs x}{k}\right) d\mu \\
&\leq \sum_{\bs p\in k\Delta \cap \bbZ^n}n_{\alpha,k}(\bs p)f(\bs p)\delta_{\alpha}\left(\dfrac{\bs x}{k}\right)
 +\dfrac{C(n+1)!}{k^2}\int_{k\Delta}a_f(\bs x)\delta_{\alpha}\left(\dfrac{\bs x}{k}\right){d\mu} \\
 & + \frac{1}{(n+2)}\sum_{\bs p\in k\Delta\cap \bbZ^n}\sum_{S\in \mathcal T \atop \bs p\in \mathcal{V}(S)}\sum_{i, j\ne i}f(\bs p)\delta_{\alpha,k, \bs v_{ij}^S}(\bs p).
\end{align*}
Summing up all $\alpha$, and using $\sum_\alpha \delta_{\alpha,k}(\bs x)=1$, 
\begin{equation}\label{eq:IntegSum}
\int_{k\Delta}a_f(\bs x)\delta_{\alpha}\left(\dfrac{\bs x}{k}\right) d\mu= \frac{1}{(n+1)!}\sum_{\bs p\in k\Delta\cap \bbZ^n} n_{\alpha,k}(\bs p)f(\bs p),
\end{equation}
we obtain the desired inequality. 
\end{proof}

The next step is to derive estimates on the boundary $\partial (k\Delta)$. However, we must keep in mind that $f$ is convex only on each face. 
%Def
\begin{definition}\rm
For a vertex $\bs p_\alpha$ of $\Delta$ and $k\in \bbZ_{>0}$, let $C(k \bs p_\alpha)$ be the vertex cone.
A \emph{simplex triangulation} $\mathcal T$ of the boundary $\dd C(k \bs p_\alpha)$ is a triangulation such that
\begin{itemize}
\item $\mathcal T$ is a simplex triangulation in the sense of Definition \ref{def:SimpTri} (2),
\item On each face $F$ of $C(k \bs p_\alpha)$, the restriction of $\mathcal T$ to $F$ is again a simplex triangulation.
\end{itemize}
We denote such a triangulation of $\dd C(k \bs p_\alpha)$ by $\mathcal T(\dd C(k \bs p_\alpha))$.
\end{definition}
Suppose that a simplex triangulation $\mathcal T(\partial C(k\bs p_\alpha))$ of $\partial (k\Delta)$ is induced from a type F triangulation of $C(k\bs p_\alpha)$. Then  it is itself a type F triangulation.

Let us extend the previous discussion to estimates on the boundary $\partial (k\Delta)$.
Define, for a vertex $\bs p_\alpha$ of $\Delta$ and $k\in \mathbb Z_{>0}$, the map 
\[
m_{\alpha,k}: (\partial C(k\bs p_\alpha)\cap \bbZ^{n-1})%\times \partial \mathcal{B}_k 
\longrightarrow \mathbb Z
\]
by
\[
m_{\alpha,k}(\bs q):= \#\set{S\in \mathcal T(\partial C(k\bs p_\alpha)) | S \text{ is an $(n-1)$-simplex with } \bs q\in \mathcal V(S)}.
\]
%$m_{\alpha,k}$ by
%\begin{equation*}
%\xymatrix@R=-1ex{
%m_{\alpha,k}: (\partial C(k\bs p_\alpha)\cap \bbZ^{n-1})\times \partial \mathcal{B}_k \ar[r]& \bbZ\\
%\quad \! \din &\din\\
 %(q,B) \ar@{|->}[r]& \#\set{S | \text{$S$ is an $(n-1)$-simplex touching $q$}}.
 %}
%\end{equation*}  
Let $R$ denote the number of the vertices of $\Delta$, and define 
\[
m_k(\bs q):=\max_{ 1\le \alpha \le R }m_{\alpha ,k}(\bs q).
\]
%Cor
Next, define the set of edge vectors in the triangulation $\mathcal T:=\mathcal T(\partial C(k\bs p_\alpha))$ by 
\[
\mathcal{V}(\mathcal T):=\Set{\bs v_{ij}=\bs p_j-\bs p_i |\bs  p_i,  \bs p_j\in S, ~~  S\in \mathcal T},
\]
i.e., all vectors defined by differences of any two vertices in a simplex $S\in \mathcal T$, taken over all simplices in the triangulation.
%Cor
\begin{corollary}\label{boundary estimates}
Let $k\in \mathbb Z_{>0}$, and let $f:k\Delta\rightarrow \bbR_{\ge 0}$ be non-negative and convex on each face of $\partial (k\Delta)$. Then there exists a constant $M_1>0$, independent of $k$, such that
\begin{align*}\label{ineq:BoundEst}
&\int_{ \dd (k\Delta)}f(\bs x){d\sigma}\\
\leq& \left(1+\frac{C}{k^2}\right)\sum_{\bs p\in \dd (k\Delta)\cap \bbZ^n}\frac{m_{k}(\bs p)f(\bs p)}{n!}+\frac{M_1}{k}\max_{\bs p\in\dd (k\Delta)\cap \bbZ^n}m_k(\bs p) \sum_{\bs p\in \dd (k\Delta)\cap \bbZ^n}f(\bs p),
\end{align*}
where $C$ is the constant defined in Corollary \ref{coro: global estimate refine}.
Moreover, the constant $M_1$ is given by
\[
M_1=\frac{R n}{(n+1)!} \max_{\bs x\in \Delta}\Norm{ \nabla \delta_{\alpha}|_{\bs x}}  C_1,
\] 
for some positive constant 
\begin{equation}\label{eq:FinBound}
C_1\geq \max_{\bs v\in\mathcal{V}(\cT)} \Norm{\bs v}.
\end{equation}
Here  $\delta_\alpha$ is the partition of unity from Section \ref{sec:partition}, and $ \nabla \delta_{\alpha}|_{\bs x}$ denotes
the gradient of $\delta_\alpha$ evaluated at the point $\bs x$.
The existence of the finite bound $C_1$ in $\eqref{eq:FinBound}$ is guaranteed by the type F triangulation.
\end{corollary}
\begin{proof} 
Let $f$ be a non-negative function whose restriction to each face of $\partial (k\Delta)$ is convex. 
Using the same argument as in Corollary \ref{coro: global estimate refine} applied to each $(n-1)$-simplex in the triangulation
of $\partial C(k\bs p_\alpha)$, we obtain
\begin{align*}
 \int_{\dd (k\Delta)}  f(\bs x){d\sigma}&
\leq  \frac{1}{n!}\sum_{\bs p\in \dd (k\Delta)\cap \bbZ^n}m_{k}(\bs p)f(\bs p)+
\frac{C}{k^2n!}\sum_{\bs p\in \dd (k\Delta)\cap \bbZ^n}m_{k}(\bs p)f(\bs p)\\
&+ \frac{1}{(n+1)!}\sum_{\alpha=1}^R\sum_{\bs p\in \dd (k\Delta)\cap \bbZ^n}\sum_{S\in \mathcal T \atop \bs p\in \mathcal{V}(S)}\sum_{i,j=1\atop j\neq i}^{n}f(\bs p)\delta_{\alpha,k, \bs v_{ij}^S}(\bs p),
\end{align*}
where $\delta_{\alpha,k, \bs v_{ij}^S}(\bs p):=\braket{\nabla \delta_{\alpha,k}(\bs p), \bs v_{ij}^S}$.
For estimating the last term, we note that 
\[
(\nabla\delta_{\alpha,k,\bs v_{ij}^S})|_{\bs p}:=\Braket{ \nabla \delta_{\alpha,k}, \bs v_{ij}^S}.
\] 
In particular, $\bs v_{ij}^S$ is a fixed vector associated with simplex $S$ and $\delta_{\alpha,k}$ is linear on each simplex, the directional derivative along $\bs v_{ij}^S$ is constant. Therefore, it does not depend on the point $\bs p$ as long as $\bs p\in \mathcal V(S)$.

By the scaling relation $\delta_{\alpha, k}(\bs x)=\delta_\alpha(\bs x/k)$, we have
\[
\nabla \delta_{\alpha,k}(\bs x)=\frac{1}{k}\nabla \delta_\alpha \left(\frac{\bs x}{k}\right)
\]
Hence, for all $\bs p\in \partial (k\Delta)$,
\begin{align*}
\norm{(\nabla\delta_{\alpha,k,\bs v_{ij}^S})|_{\bs p}}&=\norm{\braket{\nabla \delta_{\alpha,k}(\bs p), \bs v_{ij}^S}}
\leq \Norm{\nabla \delta_{\alpha,k}(\bs p)}\cdot \Norm{\bs v_{ij}^S} \\ %\max_{\bs v\in\mathcal{V}(\cT)} \Norm{\bs v}
&\leq \max_{\bs x\in k\Delta}\Norm{\nabla \delta_{\alpha,k}(\bs x)}  \max_{\bs v\in\mathcal{V}(\cT)} \Norm{\bs v}
=\dfrac{1}{k} \max_{\bs x\in \Delta}\Norm{ \nabla \delta_{\alpha}(\bs x)}  \max_{\bs v\in\mathcal{V}(\cT)} \Norm{\bs v}.
\end{align*} 
From this, as $f\geq 0$, we obtain
\begin{align*}
 \int_{\dd (k\Delta)}  f(\bs x){d\sigma}&
\leq& \frac{1}{n!} \left(1+\dfrac{C}{k^2}\right)\sum_{\bs p\in \dd k\Delta\cap \bbZ^n}m_{k}(\bs p)f(\bs p)+\dfrac{M_1}{k}\sum_{\bs p\in \dd (k\Delta)\cap \bbZ^n}f(\bs p),
\end{align*}
where 
\[
M_1=\frac{R n}{(n+1)!} \max_{\bs x\in \Delta}\Norm {\nabla \delta_{\alpha}|_{\bs x}}  \max_{\bs v\in\mathcal{V}(\cT)} \Norm{\bs v}.
\]
Type F triangulation implies that $\max_{\bs v\in\mathcal{V}(\cT)} \Norm{\bs v}\leq C_1$ for some $C_1>0$.
\end{proof}

%Sec2.5
\subsection{Ideal triangulation}\label{sec:Ideal}
We now simplify the term 
\[
\frac{1}{(n+2)!}
\sum_{\alpha=1}^R\sum_{\bs p\in k\Delta\cap \bbZ^n}\sum_{S\in \mathcal T \atop 
\bs p\in \mathcal{V}(S)}\sum_{i,j=1\atop j\neq i}^{n+1}f(\bs p)\delta_{\alpha,k, \bs v_{ij}^S}(\bs p)
\] 
appearing in Corollary \ref{coro: global estimate refine}.
%Def
\begin{definition}\rm
Let $\bs p_\alpha$ be vertex of $\Delta$ and let $C(\bs p_\alpha)$ be its vertex cone.
 A type F triangulation $\mathcal T(C(\bs p_\alpha))$ is called {\emph{strongly ideal}}, if for any interior lattice point 
 $\bs p \in (C(\bs p_\alpha))^o \cap \mathbb Z^n$, and any $S\in \mathcal T:=\mathcal T(C(\bs p_\alpha))$ containing $\bs p$ as a vertex, it satisfies
\[
\sum_{S\in \mathcal T \atop \bs p\in \mathcal{V}(S)}\sum_{\bs q\in \mathcal{V}(S)}(\bs q-\bs p)=0.
\]
We call a triangulation $\mathcal T$ of $\Delta$ is {\emph{strongly ideal}} if every vertex cone admits a strongly ideal triangulation.
Moreover, a type F triangulation $\mathcal T$ is called {\emph{ideal}} if it is a finite union of strongly ideal triangulations. 
\end{definition}
%Lem
\begin{lemma}\label{lem:StrIdeal}
Let $C(\bs p_\alpha)$ be a vertex cone of $\Delta$.
A type F triangulation $\mathcal T(C(\bs p_\alpha))$ is strongly ideal if and only if, for any interior lattice point 
$\bs p\in (C(\bs p_\alpha))^o \cap \mathbb Z^n$, any function $f$, and any smooth function $\delta_{\alpha}$ on 
$(C(\bs p_\alpha))^o$, 
\[
\sum_{S\in T \atop 
\bs p\in \mathcal{V}(S)}\sum_{\bs q\in \mathcal{V}(S)}f(\bs p)\delta_{\alpha, \bs q- \bs p}(\bs p)=0.
\]
\end{lemma}
\begin{proof}
We compute
   \begin{align*}
   \sum_{S\in \mathcal T \atop \bs p\in \mathcal{V}(S)}\sum_{\bs q\in \mathcal{V}(S)}f(\bs p)\delta_{\alpha, \bs q-\bs p}(\bs p)
    &=    \sum_{S\in \mathcal T \atop \bs p\in \mathcal{V}(S)}\sum_{\bs q\in \mathcal{V}(S)}f(\bs p)\braket{ \nabla \delta_{\alpha}(\bs p),\bs q-\bs p}\\
    =& f(\bs p)\left\langle \nabla \delta_{\alpha}(\bs p),\sum_{S\in \mathcal T \atop \bs p\in \mathcal{V}(S)}\sum_{\bs q\in \mathcal{V}(S)}(\bs q-\bs p)\right\rangle.
   \end{align*} 
The statement follows immediately from the definition of a strongly ideal triangulations. 
\end{proof}
We next introduce a special class of triangulations, which forms a subclass of strongly ideal triangulations.
%Def
\begin{definition}\rm
 A type F triangulation $\mathcal T(C(\bs p_\alpha))$ is called {\emph{locally centrally symmetric}} if, for any interior lattice point $\bs p\in   (C(\bs p_\alpha))^o\cap \mathbb Z^n$, and any simplex $S\in \mathcal T$ satisfying $\bs p\in \mathcal V(S)$, one has 
\begin{equation}\label{def:LocSym}
-(S-\bs p)+\bs p\in \mathcal T.
\end{equation}
\end{definition}
%Lem
\begin{lemma}\label{lem: vanishing 1/k term}
Every locally centrally symmetric triangulation is strongly ideal. 
\end{lemma}
\begin{proof}
Let $\mathcal T$ be a locally centrally symmetric triangulation of $\Delta$.
By Lemma \ref{lem:StrIdeal}, it suffices to prove that 
\[
\sum_{S\in \mathcal T \atop \bs p\in \mathcal{V}(S)}\sum_{\bs q\in \mathcal{V}(S)}f(\bs p)\delta_{\alpha, \bs q-\bs p}(\bs p)=0.
\]
for any interior lattice point $\bs p\in (C(\bs p_\alpha))^o \cap \mathbb Z^n$.

Fix such a point $\bs p$.
By the symmetry condition $\eqref{def:LocSym}$, for every simplex $S\in \mathcal T$ with $\bs p\in \mathcal V(S)$,
the reflected simplex
\[
S':=-(S-\bs p)+\bs p
\]
also belongs to $\mathcal T$. 
Moreover, for any vertex $\bs q\in \mathcal V(S)$, the corresponding vertex of $S'$ is
\[
\bs r:=-(\bs q -\bs p)+\bs p.
\]
Thus,
\[
\bs q\in \mathcal V(S) \quad \Leftrightarrow \quad \bs r \in \mathcal V(S').
\]
This yields that for $\bs p \in \Delta^o\cap \mathbb Z^n$, we have
\[
\sum_{\bs q\in \mathcal V(S)}f(\bs p)\delta_{\alpha, \bs q -\bs p}(\bs p)=\sum_{\bs r\in \mathcal V(S')}f(\bs p)\delta_{\alpha, \bs r -\bs p}(\bs p).
\]
Thus, 
 \begin{align*}
 \sum_{\bs q\in \mathcal{V}(S)}f(\bs p)\delta_{\alpha, \bs q-\bs p}(\bs p)
 =&\dfrac{1}{2}\left(\sum_{\bs q\in \mathcal{V}(S)}f(\bs p)\delta_{\alpha, \bs q-\bs p}(\bs p)
 +\sum_{r\in \mathcal{V}(S')}f(\bs p)\delta_{\alpha, \bs r-\bs p}(\bs p)\right)\\
=&\frac{1}{2}\left(\sum_{\bs q\in \mathcal{V}(S)}f(\bs p)\delta_{\alpha, \bs q-\bs p}(\bs p)
+\sum_{\bs q\in \mathcal{V}(S)}f(\bs p)\delta_{\alpha, -(\bs q-\bs p)}(\bs p)\right).
 \end{align*}
 Summing over all $S\in \mathcal T$ with $\bs p\in \mathcal V(S)$, we obtain
 \begin{align*}
 \sum_{S\in \mathcal T \atop \bs p\in \mathcal{V}(S)}\sum_{\bs q\in \mathcal{V}(S)}f(\bs p)\delta_{\alpha, \bs q-\bs p}(\bs p)
 =&\dfrac{1}{2}\left(\sum_{S\in \mathcal T \atop \bs p\in \mathcal{V}(S)}\sum_{\bs q\in \mathcal{V}(S)}f(\bs p)(\delta_{\alpha, \bs q-\bs p}(\bs p)+\delta_{\alpha, -(\bs q-\bs p)}(\bs p))\right)\\
=&\dfrac{1}{2}\left(\sum_{S\in \mathcal T \atop \bs p\in \mathcal{V}(S)}\sum_{\bs q\in \mathcal{V}(S)}f(\bs p)(\langle \nabla \delta_{\alpha}(\bs p),(\bs q-\bs p)-(\bs q-\bs p))\rangle\right)\\
=&0.
 \end{align*}
\end{proof}
%Remark
\begin{remark}\label{rem:lsc}\rm
Locally centrally symmetric triangulations are easy to construct.
Let 
\[
Q=\smallcup_i A_i
\]
be the union of all boundary faces of simplices in a triangulation $\mathcal T$ of a cone $C$, where each $A_i$ is a connected flat domain contained in an affine subspace. 
If for every interior lattice points $\bs p\in \Delta^o\cap \mathbb Z^n$, the condition $\bs p\in A_i$ implies $\bs p\in \mathrm{Int}(A_i)$. Then $\mathcal T$ is locally centrally symmetric. In particular, the standard triangulations in Section \ref{sec: sct} satisfy this condition and are therefore strongly ideal. 
\end{remark}
As an immediate application of Corollary \ref{coro: global estimate refine}, we obtain the following result.
%Cor
\begin{corollary}\label{cor:application}
Let $\Delta=\conv\set{\bs p_{\alpha}}_{\alpha=1}^R$ be an integral polytope with vertices $\bs p_{\alpha}$.
Let $f$ be a non-negative convex function on $k\Delta$, and let $\Set{\delta_{\alpha}}_{\alpha=1}^R$ be a partition of unity as above.
Assume that each vertex cone $C(k \bs p_{\alpha})$ of $k\Delta$ admits a strongly ideal triangulation. Then 
\begin{align*}
 &\left(1-\dfrac{C}{k^2}\right)\int_{k\Delta} f(\bs x){d\mu}\\
\leq&\frac{1}{(n+1)!}\sum_{\bs p\in k\Delta\cap \bbZ^n}n_{k}(\bs p)f(\bs p)+\frac{M}{k}
\max_{\bs p\in\partial (k\Delta)\cap \bbZ^n}n_k(\bs p)\sum_{\bs p\in \partial (k\Delta)\cap \bbZ^n}f(\bs p),
\end{align*}
where 
$n_k(\bs p)=\max_{\alpha} n_{\alpha,k}(\bs p)$, and $C$ is the constant defined in $\eqref{def:constC}$.
The constant $M>0$ is given by 
\begin{equation}\label{eq:ConstM}
M=C'\max_{\bs x\in \Delta}\Norm {\nabla \delta_{\alpha}({\bs x}) } \max_{\bs v\in\mathcal{V}(\cT)} \Norm{\bs v},
\end{equation}
where $C'>0$ depends only on the combinatorics of the triangulation, and  $\max_{\bs v\in\mathcal{V}(\cT)} \Norm{ \bs v}$ is bounded because $\cT$ is a type F triangulation.
\end{corollary}
\begin{proof}
Let $\mathcal T$ be a strongly ideal triangulation as assumed.
By Corollary \ref{coro: global estimate refine}, we have 
\begin{align}\label{ineq:StrIdeal}
\int_{k\Delta}  f(\bs x){d\mu} &\leq  \dfrac{1}{(n+1)!}\sum_{\bs p\in k\Delta\cap \bbZ^n}n_{k}(\bs p)f(\bs p)
+\dfrac{C}{k^2(n+1)!}\sum_{\bs p\in k\Delta\cap \bbZ^n}n_{\alpha,k}(\bs p)f(\bs p)\\
& +\frac{1}{(n+2)!}\sum_{\alpha=1}^R\sum_{\bs p\in k\Delta\cap \bbZ^n}\sum_{S\in \mathcal T  \atop \bs p\in \mathcal{V}(S)}
\sum_{i, j\neq i}^{n+1}f(\bs p)\delta_{\alpha,k, \bs v_{ij}^S}(\bs p).
\end{align}
Since $\mathcal T$ is strongly ideal, the interior contribution vanishes by Lemma \ref{lem:StrIdeal}:
\[
\sum_{S\in \mathcal T  \atop \bs p\in \mathcal{V}(S)}
\sum_{i, j\neq i}^{n+1}f(\bs p)\delta_{\alpha,k, \bs v_{ij}^S}(\bs p)=0 \quad \text{ for all } \quad \bs p\in (k\Delta)^o\cap \bbZ^n.
\]
Thus, the last term reduces to a sum over boundary lattice points:
\[
\sum_{\alpha=1}^R\sum_{\bs p\in \partial(k\Delta)\cap \bbZ^n}\sum_{S\in T  \atop \bs p\in \mathcal{V}(S)}
\sum_{i, j\neq i}^{n+1}f(\bs p)\delta_{\alpha,k, \bs v_{ij}^S}(\bs p).
\]
Using the scaling relation
\[
\nabla \delta_{\alpha,k}(\bs x)=\frac{1}{k}\nabla\delta_\alpha\left( {\bs x}/{k}\right),
\qquad \bs x \in \partial(k\Delta),
\]
and applying the argument as in the proof of  Corollary \ref{boundary estimates}, we obtain
\[
\left|\delta_{\alpha,k, \bs v_{ij}^S}(\bs p)\right| \leq \dfrac{M}{k} \quad \text{ for all } \bs p \in \partial(k\Delta)\cap \bbZ^n,
\]
where $M$ is the constant defined in $\eqref{eq:ConstM}$. Consequently, we have
\[
\frac{1}{(n+2)!}\sum_{\alpha=1}^R\sum_{\bs p\in k\Delta\cap \bbZ^n}\sum_{S\in \mathcal T  \atop \bs p\in \mathcal{V}(S)}
\sum_{i, j\neq i}^{n+1}f(\bs p)\delta_{\alpha,k, \bs v_{ij}^S}(\bs p)
\leq \frac{M}{k}\max_{\bs p\in\partial (k\Delta)\cap \bbZ^n}n_k(\bs p)\sum_{\bs p\in \partial (k\Delta)\cap \bbZ^n}f(\bs p).
\]
Combining this estimate with $\eqref{ineq:StrIdeal}$ and applying $\eqref{eq:IntegSum}$ to the non-negative convex function $f$,
we obtain the desired inequality. 
\end{proof}
\noindent {\bf{Geometric Idea.}} For a strongly ideal triangulation, the interior contribution vanishes, and therefore
\[
\sum_{\alpha=1}^R\sum_{\bs p\in k\Delta\cap \mathbb Z^n}\sum_{S\in \mathcal T  \atop \bs p\in \mathcal{V}(S)}
\sum_{i, j\neq i}^{n+1}f(\bs p)\delta_{\alpha,k, \bs v_{ij}^S}(\bs p)=\sum_{\alpha=1}^R\sum_{\bs p\in \partial(k\Delta)\cap \mathbb Z^n}\sum_{S\in \mathcal T  \atop \bs p\in \mathcal{V}(S)}
\sum_{i, j\neq i}^{n+1}f(\bs p)\delta_{\alpha,k, \bs v_{ij}^S}(\bs p).
\]
For an ideal triangulation, which is a finite union of strongly ideal triangulations, additional contributions arise only from lattice points contained in the ``gluing region", i.e., along interfaces between the strongly ideal pieces. Thus, it suffices to estimate the contribution from such points (e.g: the lattice points lying on red segments in Figure \ref{fig:gluing}).
\begin{figure}[h!]
\tdplotsetmaincoords{70}{20}
\begin{tikzpicture}[tdplot_main_coords]
\draw[thick,dotted] (-2,0,0) -- (0,1,0);
\draw[thick,dotted] (2,0,0) -- (0,1,0);
\draw[thick] (-2,0,0) -- (0,-1,0);
\draw[thick] (2,0,0) -- (0,-1,0);
\draw[thick] (0,-1,0) -- (0,0,1);
\draw[red,thick] (-1,0,0) -- (0,1,0);
\draw[red,thick] (0,0,0) -- (0,1,0);
\draw[red,thick] (0,1,0) -- (1,0,0);
\draw[red, thick] (2,0,0) -- (-2,0,0);
\draw[thick] (0,0,1) -- (-2,0,0);
\draw[thick] (2,0,0) -- (0,0,1);
\draw[thick,dotted] (0,1,0) -- (0,0,1);

\draw[red,thick] (-1,0,0) -- (0,-1,0);
\draw[red,thick] (0,0,0) -- (0,-1,0);
\draw[red,thick] (0,-1,0) -- (1,0,0);
\draw[red,thick] (-1,0,0) -- (0,0,1);
\draw[red,thick] (0,0,0) -- (0,0,1);
\draw[red,thick] (0,0,1) -- (1,0,0);
\end{tikzpicture}
\begin{tikzpicture}[tdplot_main_coords]
\draw[thick,dotted] (2,0,0) -- (0,1,0);
\draw[thick] (2,0,0) -- (0,-1,0);
\draw[thick] (2,0,0) -- (0,0,1);
\draw[thick,dotted] (0,1,0) -- (0,0,1);
\draw[thick] (0,-1,0) -- (0,0,1);
\draw[thick] (2,0,0) -- (0,0,-1);
\draw[thick, dotted] (0,1,0) -- (0,0,-1);
\draw[thick] (0,-1,0) -- (0,0,-1);
\draw[thick,red] (0,-1,0) -- (1,0,0);
\draw[thick,red] (0,1,0) -- (1,0,0);
\draw[thick,red] (0,0,1) -- (1,0,0);
\draw[thick,red] (0,0,-1) -- (1,0,0);
\draw[thick,red] (2,0,0) -- (1,0,0);
\end{tikzpicture}
\caption{ The convex hulls $R(0,0,1)$ and $R(2,0,0)$, each subdivided into four simplices $S_1, \dots ,S_4$ with $(1,0,0)\in S_i$ for $i=1,\dots ,4$.}\label{fig:gluing}
\end{figure}

The key observation is that every lattice point $\bs p$ in the gluing region can be projected onto the boundary along a suitable direction. More precisely, there exists a vector $\bs v$, transverse to the faces of the cone, such that the line
\[
\Set{\bs p+t\bs v| t\in \mathbb R}
\]
intersects the boundary at two points $q_1(\bs p)$ and $q_2(\bs p)$. By the convexity and non-negativity of $f$,
we obtain
\[
f(\bs p)\leq \lambda f(q_1(\bs p))+(1-\lambda)f(q_2(\bs p))
\]
for some $\lambda \in [0,1]$, and hence
\[
f(\bs p)\leq  f(q_1(\bs p))+f(q_2(\bs p)).
\]
Thus, the values of $f$ at interior lattice points in the gluing region can be controlled by its values at boundary lattice points.

Moreover, the type F condition implies a uniform bound on the number of such projections landing in any given boundary simplex. 
%simplex on the boundary. 
%Indeed, after projection along $\bs v$, any distinct lattice points remain separated, 
Distinct lattice points remain distinct under projection, so each boundary simplex contains only finitely many projected points, with a bound depending only on the combinatorics of the triangulation.
%. A uniform bound depending only on the combinatorics of the triangulation. 

Consequently, the total contribution from the gluing region is bounded by a constant multiple of
\[
\sum_{\bs p\in \partial (k\Delta)\cap \mathbb Z^n}f(\bs p),
\] 
which explains the boundary term appearing in Corollary \ref{cor:application}.

%Ex
\begin{example}\rm
In Figure \ref{fig:NonSmall}, let $f$ be a non-negative convex function. Then, for any positive integer $n$,
\[
f(-n,0)\leq \frac{1}{2}f\left(-n, \dfrac{n}{3}\right)+ \frac{1}{2}f\left(-n, -\dfrac{n}{3}\right).
\]
In particular, the value of $f(-n,0)$ can be estimated in terms of its values at the boundary points $(-n, \frac{n}{3})$ and $(-n, -\frac{n}{3})$.
\begin{figure}[h!]
\begin{tikzpicture}[x=.5cm,y=.5cm]
\draw[thick] (9,0) -- (0,-3);
\draw[thick] (9,0) -- (0,3);
\filldraw[green] (6,0) circle (2pt);
\filldraw[green] (5,0) circle (2pt);
\filldraw[green] (4,0) circle (2pt);
\filldraw[green] (3,0) circle (2pt);
\draw[yellow,thick] (3,2) -- (6,1);
\draw[thick] (3,2) -- (3,-2);
\draw[thick] (4,1.66) -- (4,-1.66);
\draw[thick] (5,1.33) -- (5,-1.33);
\draw[thick] (6,1) -- (6,-1);
\filldraw[black] (6,1) circle (2pt)node[above,right]{(-3,1)};
\filldraw[black] (6,-1) circle (2pt);
\filldraw[black] (3,2) circle (2pt)node[above]{(-6,2)};
\filldraw[black] (3,-2) circle (2pt);
\filldraw[red] (5,1.33) circle (2pt);
\filldraw[red] (5,-1.33) circle (2pt);
\filldraw[red] (4,1.66) circle (2pt);
\filldraw[red] (4,-1.66) circle (2pt);
\filldraw[black] (9,0) circle (2pt)node[anchor=west]{(0,0)};
\end{tikzpicture}
\caption{Non-small triangulation}\label{fig:NonSmall}
\end{figure}

Moreover, when projecting lattice points onto the boundary, each simplex -- namely, each segment connecting points such as 
\[
\left(-n,\frac{n}{3}\right) \quad \text{and} \quad \left(-n-1,\frac{n+1}{3}\right)
\]
-- contains at most four projected points (see the yellow highlighted segment in Figure \ref{fig:NonSmall}). In general, this number is uniformly bounded.
%(for example, the yellow segment.) And all my argument below is to show this  in general is finite. And again, for any red point, for example, \[f(-4,\frac{4}{3})\leq \dfrac{1}{3}f(-6,2)+\dfrac{2}{3}f(-3,1)\leq f(-6,2)+f(-3,1).\]

As another example, consider a lattice point in the interior (e.g., a ``red point" in the figure). By convexity of $f$,
\[
f\left(-4,\frac{4}{3}\right)\leq \dfrac{1}{3}f(-6,2)+\dfrac{2}{3}f(-3,1)\leq f(-6,2)+f(-3,1).
\]
Thus, even for interior points, the value of $f$ can be controlled by its values at boundary lattice points.

Consequently, one can estimate contributions from interior lattice points using only boundary values, which is consistent with the general strategy described below (see, the proof of Corollary \ref{coro: global estimate refine 2}).
%So we can still apply the estimate by the integral points on the boundary.
\end{example}

With an ideal triangulation, we obtain a stronger estimate. The following is the main result of this section.
%Cor
\begin{corollary}\label{coro: global estimate refine 2}
Let $\Delta=\conv\set{\bs p_{\alpha}}_{\alpha=1}^R$ be an integral polytope with vertices $\bs p_{\alpha}$.
Let $f$ be a non-negative convex function on $\Delta$, and let $\Set{\delta_{\alpha}}$ be a partition of unity as before.
Assume that each vertex cone $C(k \bs p_{\alpha})$ of $k\Delta$ admits an ideal triangulation. Then 
\begin{align*}
 &\left(1-\dfrac{C}{k^2}\right)\int_{k\Delta} f(\bs x){d\mu}\\
\leq&\frac{1}{(n+1)!}\sum_{\bs p\in k\Delta\cap \bbZ^n}n_{k}(\bs p)f(\bs p)+\frac{M}{k}
\max_{\bs p\in\dd (k\Delta)\cap \bbZ^n}n_k(\bs p)\sum_{\bs p\in \dd (k\Delta)\cap \bbZ^n}f(\bs p),
\end{align*}
where %$\delta_{\alpha,k}(\bs x)$, 
$n_k(\bs p)$ and $C$ are as in Corollary \ref{coro: global estimate refine}. 
Moreover, the constant $M>0$ is given by
\[
 M:=C'M_0rN,
%M=C'\max_{1\leq \alpha \leq R \atop\bs p \in \partial(k\Delta) \cap \mathbb Z^n}(|\nabla \delta_\alpha (\bs p)|\cdot \|\bs v_{ij}^S\|),
\]
with $C', M_0, r, N$s defined in the proof below.
%\begin{itemize}
%\item $\delta_{\alpha,k}(\bs x):=\delta_{\alpha}\left(\dfrac{\bs x}{k}\right)$,
%\item $\displaystyle n_k(\bs p)=\max_{1\leq \alpha \leq R}n_{\alpha,k}(\bs p)$,
%\item 
%\[
%C=\max_{S\in T(C(k\bs p_{\alpha}))}C_S\cdot \max_{\bs x\in \Delta}\|\mathrm{Hess}(\delta_{\alpha}(\bs x))\|,
%\]
%\delta_{\alpha,k}(\bs x)=\delta_{\alpha}\left(\dfrac{\bs x}{k}\right),
%\] 
%\item and  $M> 0$ is a positive constant.
%\end{itemize}
\end{corollary}
\begin{proof}
Fix $k\in \mathbb Z_{>0}$. For each vertex $\bs p_\alpha$, decompose the vertex cone %$C(k\bs p_\alpha)$ as 
\[
C(k\bs p_\alpha)=\smallcup_{\gamma =1}^{r_\alpha} C_\gamma,
\] 
where each $C_{\gamma}$ admits a strongly ideal triangulation. Set 
\[
r:=\max_{1\leq \alpha \leq R}r_\alpha.
\]
Let $\Set{Q_{\gamma}^{\eta_\gamma}}$ denote the $(n-1)$-dimensional faces of $C_{\gamma}$,
so that 
\[
\partial C_\gamma=\smallcup_{\eta_\gamma}Q_\gamma^{\eta_\gamma}.
\]
In particular, the collection $\Set{Q_{\gamma}^{\eta_\gamma}}$ is finite.

Fix $\gamma$ and a face $Q_\gamma^\eta$ for $\eta_\gamma=:\eta$. %We abbreviate $\eta_\gamma$ by $\eta$.
By the estimate in the strongly ideal case, there exists a constant $M_0$ such that
%we obtain
\begin{align}\label{ineq:StrIdealTri}
\sum_{\alpha=1}^R \sum_{\bs p\in Q_{\gamma}^{\eta}\cap k\Delta\cap \mathbb Z^n}
\sum_{S\in \mathcal T \atop \bs p\in \mathcal{V}(S)}\sum_{i, j\ne i}^{n+1}f(\bs p)\delta_{\alpha,k, \bs v_{ij}^S}(\bs p)
\leq 
\dfrac{M_0}{k}\max_{\bs p\in Q_{\gamma}^{\eta}\cap k\Delta\cap \mathbb Z^n} n_k(\bs p)
\sum_{\bs p\in Q_{\gamma}^{\eta}\cap k\Delta\cap \mathbb Z^n} f(\bs p).
%\sum_{S\in T \atop\bs p\in \mathcal{V}(S)}f(\bs p),
\end{align}
%for a constant $M_0$.
%with the constant $M$ in $\eqref{eq:ConstM2}$.

Next, choose a primitive vector $\bs v\in \mathbb Z^n$ not tangent to any face $Q_\gamma^\eta$.
For any lattice point $\bs p\in Q_\gamma^\eta$ not on the boundary of $Q_\gamma^\eta$, the line
\[
\Set{\bs p+t \bs v|t\in \mathbb R}
\]
intersects $\partial C_\gamma$ at two points $q_1(\bs p)$ and $q_2(\bs p)$.
By convexity and non-negativity of $f$,
\[
f(\bs p)\leq f(q_1(\bs p))+f(q_2(\bs p)).
\]
Let $\mathcal T(C_{\gamma}):=\restrict{\mathcal T}{C_\gamma}$ be the triangulation of $C_\gamma$ induced by $\mathcal T:=\mathcal T(C(k\bs p_\alpha))$.
By using $\mathcal T(C_\gamma)$, each $q_i(\bs p)$ $(i=1,2)$ lies in a simplex $S_{q_i(\bs p)}$ with vertices
$P_{ q_i(\bs p),0}, \cdots, P_{ q_i(\bs p),n-1}$, and hence
\[
f( q_i(\bs p))\leq \sum_{j=0}^{n-1}f\left(P_{ q_i(\bs p),j} \right ).
\]
It remains to %show the boundedness of 
bound the number of lattice points projecting to a given boundary simplex.
Let $\mathrm{proj}_{\bs v}$ denote projection along $\bs v$, so that $\mathrm{proj}_{\bs v}(\bs p)\in \set{q_1(\bs p), q_2(\bs p)}$. 
%Under $\mathrm{proj}_{\bs v}$, 
Using this projection, distances are uniformly controlled:
\[
d(q_i(\bs p_1), q_i(\bs p_2))\geq c \cdot d(\bs p_1, \bs p_2),
\]
for some uniform constant $c>0$. Consequently, distinct lattice points remain distinct after projection, and each boundary
simplex contains at most $C'$ projected points.
%a bounded number $C'$ of projected points. %$\mathrm{proj}_{\bs v}(\bs p_j)$.
Combining this with $\eqref{ineq:StrIdealTri}$, we obtain
\begin{align}
\begin{split}\label{ineq:bounded}
\sum_{\alpha=1}^R\sum_{\bs p\in Q_{\gamma}^{\eta}\cap k\Delta\cap \mathbb Z^n}
\sum_{S\in \mathcal T \atop \bs p\in \mathcal{V}(S)}\sum_{i, j\neq i}^{n+1}f(\bs p)\delta_{\alpha,k, \bs v_{ij}^S}(\bs p)
\leq&\frac{C'M_0}{k}\cdot rN\max_{\bs p\in\partial (k\Delta)\cap \bbZ^n}n_k(\bs p)
\sum_{p\in \dd k\Delta\cap \bbZ^n}f(\bs p),
\end{split}
\end{align}
where  
\[
N:=\max_{1\leq \gamma \leq r}\#\Set{\text{facets of }C_\gamma}.
%\mathcal F(C_\gamma)=\max_{1\leq \gamma \leq r} \eta_\gamma.
\]
Finally, combining Corollary \ref{coro: global estimate refine}, $\eqref{eq:IntegSum}$, $\eqref{ineq:bounded}$ and the scaling relation
\[
\nabla \delta_{\alpha,k}(\bs x)=\frac{1}{k}\nabla \delta_\alpha(\bs x/k),
\]
we obtain 
\begin{align*}
 &\frac{1}{(n+1)!}\sum_{\bs p\in k\Delta\cap \mathbb Z^n} n_{k}(\bs p)f(\bs p)\\
\geq&\int_{k\Delta} f(\bs x)d\mu- \frac{C}{k^2}\sum_{\alpha=1}^R\sum_{\bs p\in k\Delta \cap \mathbb Z^n}
n_{\alpha, k}(\bs p)f(\bs p)-\frac{1}{(n+2)!}\sum_{\alpha=1}^R\sum_{\bs p\in \partial (k\Delta)\cap \mathbb Z^n}
\sum_{S\in T\atop \bs p\in \mathcal V(S)}\sum_{i, j\neq i}^{n+1}f(\bs p)\delta_{\alpha, k, \bs v_{ij}^S}(\bs p)\\
\geq&\left(1-\dfrac{C}{k^2}\right)\int_{k\Delta} f(\bs x){d\mu}-\frac{C'M_0}{k}\cdot rN\max_{\bs p\in\partial (k\Delta)\cap \bbZ^n}n_k(\bs p)
\sum_{\bs p\in \dd k\Delta\cap \bbZ^n}f(\bs p),
\end{align*}
which yields the desired inequality. 
Here we replaced the constants $C$ and $M_0$ by $C/(n+1)!$  and $M_0/(n+2)!$, respectively.
\end{proof}

%Sec3
\section{Small polytopes and uniform K-stability of toric varieties}\label{sec:SmPolyUnifK}
%Sec3.1
\subsection{Small polytopes}\label{sec:SmallPolytope}
%Def
\begin{definition}\label{def:SmallPolytope}\rm
Let $\Delta$ be an integral polytope with vertex set $\mathcal V(\Delta)=\set{\bs p_1, \dots , \bs p_R}$
and let $\Delta^o$ denote its interior. Fix a positive integer $k$.

We say that the cone $C(\bs p_i)$ is {\emph{small}} with respect to a type F triangulation of $C(\bs p_i)$ and its boundary $\dd C(\bs p_i)$ if the following conditions hold:
\begin{enumerate}
\item For any interior lattice point $\bs p\in k\Delta^o \cap \bbZ^n$, %we have
\[
n_{i,k}(\bs p)\leq (n+1)!;
\] 
\item For any boundary lattice point $\bs p\in (\dd (k\Delta)\setminus \mathcal V(k\Delta))\cap \bbZ^{n}$, 
%except for the vertices, we have 
\[
n_{i,k}(\bs p)\leq \dfrac{(n+1)!}{2};
\]
\item For any boundary lattice point $\bs p\in (\dd (k\Delta)\setminus  \mathcal V(k\Delta))\cap \bbZ^{n}$ 
%except for the vertices, we have
\[
m_{i,k}(\bs p)\leq n!.
\]
\end{enumerate}
We say that $\Delta$ is a {\emph{small polytope}} if there exists a type F triangulation of each cone $C(\bs p_i)$ such that every cone $C(\bs p_i)$ is small.
\end{definition}  
From the definitions of $n_k$ and $m_k$, we immediately obtain the following. 
%Lem
\begin{lemma}\label{lem: n_k}
Let $\Delta$  be a small polytope. Then:
\begin{enumerate}
\item For any interior lattice point $\bs p\in k\Delta^o \cap \bbZ^n$, %we have
\[
n_{k}(\bs p)\leq (n+1)!;
\] 
\item For any boundary lattice point $\bs p\in (\dd (k\Delta)\setminus \mathcal V(k\Delta))\cap \bbZ^{n}$ %except for the vertices, we have 
\[
n_{k}(\bs p)\leq \dfrac{(n+1)!}{2};
\]
\item For any boundary lattice point $\bs p\in (\dd (k\Delta)\setminus \mathcal V(k\Delta))\cap \bbZ^{n}$ %except for the vertices, we have
\[
m_{k}(\bs p)\leq n!.
\]
\end{enumerate}
Here,  
\[
 n_k(\bs p)=\max_{i=1,\dots,R}n_{i,k}(\bs p), \qquad  m_k(\bs p)=\max_{i=1, \dots ,R}m_{i,k}(\bs p).
\]
\end{lemma}
\begin{proof}The statements follows immediately from the definitions of $n_k(\bs p)$, $m_k(\bs p)$, and the notion of a small polytope.
%By the definition, $\displaystyle n_k(p)=\max_{i=1,\dots ,R}n_{i,k}(p)$ and $\displaystyle m_k(p)=\max_{i=1, \dots ,R}m_{i,k}(p)$. Results follow from the definition of a small polytope.
\end{proof}
%Ex
\begin{example}\label{ex:SmallPoly}\rm
Let $\Delta$ be the integral Delzant polytope corresponding to an $n$-dimensional polarized toric manifold $(X,L)$. 
Then each vertex cone $C(k\bs p_i)$ is isomorphic to the standard cone, %which is defined as 
\[
\Pi_n:=\Set{\bs x=(x_1,\cdots, x_n)\in \bbR^n | x_i\geq 0, ~~  i=1,\dots,n}.
\]
The standard triangulation of the $k$-simplex described in Lemma \ref{lem:Lee25} induces a type F triangulation of each vertex cone and its boundary.
By Lemma \ref{standard cone is small}, every such cone is small. Consequently, $\Delta$ is a small polytope.    
\end{example}

%Ex2
\begin{example}\label{X_2}\rm
Consider the polytope
\[
\Delta:=\conv\Set{(-3,0),(3,0),(0,-1),(0,1)}.
\]
For the cone $C((3,0))$, the lattice point $(0,1)$ is adjacent to five triangles in the triangulation illustrated in Figure \ref{fig:T1}.
More generally, for any positive integer $k$, the cone $C((3k,0))$ is not small with respect to any triangulation induced by this pattern.
%%%%Figure
\begin{figure}[h!]
\begin{tikzpicture}[x=.5cm,y=.5cm]
\draw[thick] (9,0) -- (0,-3);
\draw[thick] (9,0) -- (0,3);
\draw[red,thick] (6,1) -- (6,0);
\draw[red,thick] (5,1) -- (6,0);
\draw[red,thick] (6,1) -- (7,0);
\draw[red,thick] (6,1) -- (8,0);
\draw[red,thick] (6,1) -- (0,1);
\draw[red,thick] (9,0) -- (0,0);
\node[text width=1cm] at (6.5,1.8) {$(0,1)$};
\node[text width=1cm] at (10.5,0) {$(3,0)$};
\end{tikzpicture}
\caption{Non-small triangulation}\label{fig:T1}
\end{figure}

However, if one instead triangulates the cone as in Figure \ref{fig:T2}, and applies the same construction to the remaining cones, then each cone becomes small. Hence $\Delta$ is a small polytope with respect to this modified triangulation.

\begin{figure}[h!]
\begin{tikzpicture}[x=.5cm,y=.5cm]
\draw[thick] (0,0) -- (-6,2);
\draw[thick] (0,0) -- (-6,-2);
\draw[red,thick] (0,0) -- (-1,0);
\draw[red,thick] (-1,0) -- (-3,1);
\draw[red,thick] (-3,1) -- (-4,1);
\draw[red,thick] (-4,1) -- (-1,0);
\draw[red,thick] (-4,1) -- (-6,2);
\draw[red,thick] (-4,1) -- (-7,2);
\draw[red,thick] (0,0) -- (-1,0);
\draw[red,thick] (-1,0) -- (-3,-1);
\draw[red,thick] (-3,-1) -- (-4,-1);
\draw[red,thick] (-4,-1) -- (-1,0);
\draw[red,thick] (-4,-1) -- (-6,-2);
\draw[red,thick] (-4,-1) -- (-7,-2);
\end{tikzpicture}
\caption{Small triangulation}\label{fig:T2}
\end{figure}
\end{example}

%Sec3.2
\subsection{Uniform K-stability of toric varieties}\label{sec:UniformK}
Let $\bs p$ be a (not necessarily lattice) point in an integral polytope $\Delta$.
%Recall that a normalized convex function (at $p$) is a convex function $f(x)$ such that   
A convex function $f:\Delta \to \mathbb R$ is said to be {\emph{normalized at}} $\bs p$ if
\[
f(\bs p)=\min_{\bs x\in \Delta}f(\bs x)=0.
\] 
Let $O$ denote the barycenter of $\Delta$. Since translating $\Delta$ does not affect uniform K-stability,
we may assume without loss of generality that $\bs p=O$. 

%A polarized toric variety is uniformly K-stable iff the following holds: 
The following criterion characterizes uniform K-stability for polarized toric varieties.
%Def
\begin{definition}\label{def:Lambda_Sta} \rm (\cite{CLS14}; see also Proposition 3.6 in \cite{His20}). 
Let 
\[
a:=\frac{\vol(\dd \Delta, d\sigma)}{\vol(\Delta)},
\]
and let $f:\Delta \rightarrow \bbR$ be a convex function. Define
\[
\mathcal L_a(f):=\int_{\dd \Delta}f(\bs x)d\sigma-a\int_{\Delta}f(\bs x){d\mu}.
\] 
%Let $\lambda$ be a real number. Then, 
We say that $\Delta$ is {\emph{uniformly K-stable}} if the following conditions hold:
\begin{enumerate}
\item[(i)] $\mathcal L_a(\ell)=0$ for every affine function $\ell(\bs x)=a_1x_1+\dots +a_nx_n+c$; 
\item[(ii)] There exists a constant $\lambda>0$ such that, for every non-affine convex function $f$ normalized at $O$, %we have the inequality 
\begin{equation}\label{ineq:lambda}
\mathcal L_a(f)\geq \lambda \int_{\dd \Delta}f(\bs x)d\sigma.
\end{equation}
\end{enumerate}
We note that $\eqref{ineq:lambda}$ is equivalent to 
\[
a\int_{\Delta}f(\bs x){d\mu}\leq (1-\lambda)\int_{\dd \Delta}f(\bs x)d\sigma.
\]
\end{definition}

%Sec4
\section{Combinatorial criteria for asymptotic Chow polystablity}\label{sec:Criteria}
In this section, we provide a sufficient combinatorial criterion for a toric variety corresponding to a polytope $\Delta$ to be asymptotically Chow polystable.

%Sec4.1
\subsection{Main Statements}
By \cite{Ono13} (see also \cite{LLSW19}), we recall the following criterion for asymptotic Chow stability of a polarized toric variety.  
This also yields an obstruction to asymptotic Chow semistability of a polarized toric variety $(X,L)$, as observed in \cite{Ono11}.
%Def
\begin{definition}\label{eq:Chowss} \rm
Let $(X,L)$ be a polarized toric variety with the moment polytope $\Delta$. 
For a positive integer $k$ and an affine function $\ell$, we define the {\emph{Futaki-Ono invariant}} by
%We define the Futaki Ono invariant with respect to $(\Delta, L, k)$ for the affine function $l$ is given by 
\begin{align*}
FO_k(\ell):=&\sum_{\bs p
\in k\Delta\cap \mathbb{Z}^n}\ell (\bs p)-\dfrac{\chi(k\Delta)}{\vol({k\Delta})}\int_{k\Delta}\ell(\bs x)d\mu,
\end{align*}
where $\chi(k\Delta)=\#(k\Delta\cap \bbZ^n)$ denotes the number of lattice points in $k\Delta$, and $\vol(k\Delta)=\int_{k\Delta} d\mu$.
% is the volume of $k\Delta$.
 \end{definition}

Let $\bs s_\Delta(t)$ be the sum polynomial defined in $\eqref{eq:SumPoly}$, and let $E_\Delta(t)$ be the Ehrhart polynomial of $\Delta$.
For a fixed integer $k$, one verifies that the condition
\[
FO_k(\ell)=0
\]
holds for all affine functions $\ell(\bs x)$ if and only if
%we observe that the equality $\eqref{eq:Chowss}$ holds for any affine function $\ell(x)$ if and only if
\begin{equation}\label{eq:ChowWt}
\vol(\Delta)\bs s_{\Delta}(k)-kE_{\Delta}(k)\int_\Delta \bs x\, {d\mu}=0.
\end{equation} 
%Thus, $\eqref{eq:Chowss}$ holds for affine functions and all $k$ if and only if $\eqref{eq:ChowWt}$ holds for all $k$.
Expanding $\eqref{eq:ChowWt}$, we obtain
\[
\frac{k^n}{2}\left(\vol(\Delta)\int_{\dd \Delta} \bs x\, {d\mu}-\vol(\dd \Delta, \sigma)\int_{\Delta} \bs x\, {d\mu}\right)+\sum_{j=1}^{n-1}k^j\left( \bs s_{\Delta,j}\vol(\Delta)-E_{\Delta,j-1}\int_{\Delta}\bs x\, {d\mu}\right)=0.
\]
Define
\[
\mathcal F_{\Delta,j}:=\vol(\Delta)\bs s_{\Delta,j}\vol(\Delta)-E_{\Delta,j-1}\int_{\Delta}\bs x\, {d\mu}, \qquad j=1, \dots, n.
\]
We have the following.
%Prop
\begin{proposition}\label{prop:FO_inv}
Let $(X,L)$ be a K-semistable polarized toric variety with its moment polytope $\Delta$. Then the following conditions are equivalent:
\begin{enumerate}
\item For all $k\in \bbN$ and all linear functions $\ell(\bs x)$, $\eqref{eq:Chowss}$ holds;
\item  For all $k\in \bbN$ and all affine functions $\ell(\bs x)$, $\eqref{eq:Chowss}$ holds;
\item $\mathcal F_{\Delta,j}=0$ for all $j=1,2, \dots, n$.
\end{enumerate}
\end{proposition}
\begin{proof}
(1) $\Leftrightarrow$ (2). This follows since constant functions do not contribute: for any constant $c$,  
\[
\frac{1}{\#(k\Delta\cap \bbZ^n)}\sum_{\bs p\in k\Delta \cap \bbZ^n} c-\frac{1}{\vol(k\Delta)}\int_{k\Delta}c\, {d\mu}=c-c=0.
\]

\noindent (1) $\Leftrightarrow$ (3). From the above expansion, $\eqref{eq:ChowWt}$ is equivalent to
\[
\frac{k^n}{2}\left(\vol(\Delta)\int_{\dd \Delta} \bs x\, {d\mu}-\vol(\dd \Delta, \sigma)\int_{\Delta} \bs x\, {d\mu}\right)+\sum_{j=1}^{n-1}\mathcal F_{\Delta,j}k^j=0.
\]
Since $(X,L)$ is K-semistable, we have
\[
\vol(\Delta)\int_{\dd \Delta} \bs x\,{d\mu}=\vol(\dd \Delta, \sigma)\int_{\Delta} \bs x\, {d\mu}.
\] 
%Thus, we conclude that $(1)$ is equivalent to $(3)$. 
Hence, the leading term vanishes, and the above identity holds for all $k$ if and only if
\[
\mathcal F_{\Delta,j}=0 \quad \text{for all } j=1, \dots, n.
\] 
\end{proof} 
Motivated by Proposition \ref{prop:FO_inv}, for each positive integer $k$, we define %{\emph{the Futaki-Ono invariant $FO(\ell;k)$}} by
\begin{equation}\label{def:FOInv}
FO(\ell;k):=\frac{1}{ \#(k\Delta\cap \bbZ^n)}\sum_{\bs p\in k\Delta \cap \bbZ^n}\ell(\bs p)-\frac{1}{\vol(k\Delta)}\int_{k\Delta}\ell(\bs x){d\mu},
\end{equation}
for any affine function $\ell(\bs x)$.

As discussed in $\eqref{eq:LinHull}$, the family of invariants $\mathcal F_{\mathrm{Td}^p}$ provides an obstruction (i.e., a necessary condition) for a polarized variety $(X,L)$ to be asymptotically Chow semistable, as shown by Futaki \cite{Fut04}.
In the toric setting, the formulation $\eqref{def:FOInv}$ was introduced by Ono \cite{Ono13}.
%Rem
\begin{remark}\rm
Since
\[
\vol(k\Delta)\sum_{\bs p\in k\Delta\cap \bbZ^n}\ell(\bs p)-\#(k\Delta\cap \bbZ^n)\int_{k\Delta}\ell(\bs x){d\mu}
\]
can be regarded as an $\bbR^n$-valued polynomial in $k$ of degree at most $n$.
In particular, for each coordinate function $\ell(\bs x)=x_i$ $(1\leq i\leq n)$, the function $FO(\ell;k)$ has at most $n$-roots.
%$\alpha_1^{(i)}, \dots, \alpha_n^{(i)}$ for each $\ell(x)=x_i$ with $1\leq i \leq n$.
Therefore, if $FO(\ell;k)=0$ for $n+1$ distinct integers $k=k_0, \dots, k_n$, then it vanishies identically in $k$.

Moreover, if all the Donaldson-Futaki invariants vanish, then $FO(\ell;k)$ becomes an $\bbR^n$-valued polynomial of degree at most $n-1$.
In this case, it suffices to verify the vanishing for $n$ distinct integers $k=k_1, \dots ,k_n$. %instead of $n+1$ distinct integers.
\end{remark}
The main theorem of this paper is the following.
%Thm
\begin{theorem}\label{theo: main theorem 1}
Let $(X,L)$ be a polarized toric variety with the associated integral moment polytope $\Delta$. Suppose that: 
\begin{enumerate}
\item $(X,L)$ is uniformly K-stable;
\item all Futaki-Ono invariants vanish;
\item for each vertex $\bs p_\alpha \in \mathcal V(\Delta)$, the vertex cone $C(\bs p_{\alpha})$ and its boundary $\partial C(\bs p_{\alpha})$ admit an ideal triangulation $\mathcal T$ and a type F triangulation $\partial \mathcal T:=\restrict{\mathcal T}{\partial C(\bs p_{\alpha})}$, respectively, such that $C(\bs p_{\alpha})$ is small with respect to $\mathcal T$ and $\partial \mathcal T$.
\end{enumerate}
Then $(X,L)$ is asymptotically Chow polystable.
\end{theorem}
\begin{proof}
 First, let $\set{\bs p_{\alpha}}_{\alpha=1}^R$ be the set of vertices of $\Delta$, and let $\set{\delta_{\alpha}}_{\alpha=1}^R$ be a partition of unity on $\Delta$ such that $\delta_{\alpha}(\bs p_{\beta})=\delta_{\beta}^{\alpha}$. Let $k\in \mathbb Z_{>0}$ be a positive integer.
Define  
 \[
 \delta_{\alpha,k}(\bs x):=\delta_{\alpha}\left(\dfrac{\bs x}{k}\right), \qquad \bs x\in k\Delta.
 \] 
Let $f:k\Delta\rightarrow \bbR$ be a convex function, and denote $\chi(k\Delta):=\#(k\Delta \cap \mathbb Z^n)$.
We compare the discrete and continuous averages: 
\begin{equation}\label{eq:Diff}
\dfrac{1}{\chi(k\Delta)}\sum_{\bs p\in k\Delta \cap \bbZ^n} f(\bs p)-\dfrac{1}{\vol(k\Delta)}\int_{k\Delta}f(\bs x){d\mu}.
\end{equation}
Let $C$ and $M$ be the constants in Corollary \ref{coro: global estimate refine 2}. Then
 \begin{align*}
\eqref{eq:Diff}=&\dfrac{1}{\chi(k\Delta)}\sum_{\bs p\in k\Delta \cap \bbZ^n}f(\bs p)-\dfrac{1}{\chi(k\Delta)}\left(1-\dfrac{C}{k^2}\right)\int_{k\Delta}f(\bs x){d\mu}\\
  &+\dfrac{1}{\chi(k\Delta)}\left(1-\dfrac{C}{k^2}\right)\int_{k\Delta}f(\bs x){d\mu}-\dfrac{1}{\vol(k\Delta)}\int_{k\Delta}f(\bs x){d\mu}\\
  \geq& \frac{1}{\chi(k\Delta)}\sum_{\bs p\in k\Delta \cap \bbZ^n}f(\bs p)-\frac{1}{\chi(k\Delta)}\cdot \frac{1}{(n+1)!}\sum_{\bs p\in k\Delta\cap \bbZ^n}n_{k}(\bs p)f(\bs p)\\
     &- \dfrac{1}{\chi(k\Delta)}\frac{M}{k}\max_{\bs p\in\dd (k\Delta)\cap \bbZ^n}n_k(\bs p)\sum_{\bs p\in \dd (k\Delta)\cap \bbZ^n}f(\bs p)\\
&-\dfrac{C}{k^2\chi(k\Delta)}\int_{k\Delta}f(\bs x){d\mu} +\dfrac{1}{\chi(k\Delta)}\int_{k\Delta}f(\bs x){d\mu}-\dfrac{1}{\vol(k\Delta)}\int_{k\Delta}f(\bs x){d\mu}.
 \end{align*}
Since all cones $C(\bs p_\alpha)$ are small, Lemma \ref{lem: n_k} gives:
\begin{itemize}
\item $n_k(\bs p)\leq (n+1)!$ for interior lattice points,
\item $n_k(\bs p)\leq \frac{(n+1)!}{2}$ for boundary lattice points.
\end{itemize}
Set $M':=M\cdot \max_{\bs p\in \partial(k\Delta)\cap \mathbb Z^n}n_k(\bs p)$. Then we obtain a lower bounda for
$\eqref{eq:Diff}$:
\begin{align*}
 \eqref{eq:Diff} \geq&\dfrac{1}{\chi(k\Delta)}\sum_{\bs p\in k\Delta \cap \bbZ^n}f(\bs p)-\dfrac{1}{\chi(k\Delta)}\sum_{\bs p\in k\Delta^o\cap \bbZ^n}f(\bs p)-\frac{1}{2\chi(k\Delta)}\sum_{\bs p\in \dd(k\Delta)\cap \bbZ^n}f(\bs p)\\
     &- \frac{1}{\chi(k\Delta)}\cdot\frac{M'}{k}\sum_{\bs p\in \dd( k\Delta)\cap \bbZ^n}f(\bs p)
     -\dfrac{C}{k^2\chi(k\Delta)}\int_{k\Delta}f(\bs x){d\mu} \\
 &+\dfrac{1}{\chi(k\Delta)}\int_{k\Delta}f(\bs x){d\mu}-\dfrac{1}{\vol(k\Delta)}\int_{k\Delta}f(\bs x){d\mu} \\
     =&\dfrac{1}{2\chi(k\Delta)}\sum_{\bs p\in \dd(k\Delta)\cap \bbZ^n}f(\bs p)- \dfrac{M'}{k\chi(k\Delta)}\sum_{\bs p\in \dd(k\Delta)\cap \bbZ^n}f(\bs p)-\dfrac{C}{k^2\chi(k\Delta)}\int_{k\Delta}f(\bs x){d\mu}\\
     &+\left( \frac{1}{\chi(k\Delta)} -\frac{1}{\vol(k\Delta)} \right) \int_{k\Delta}f(\bs x) d\mu.
 \end{align*}
On the other hand, by the Ehrhart expansion,  
 \[
 \chi(k\Delta)=k^n\vol(k\Delta)+\frac{k^{n-1}}{2}\vol(\dd \Delta, d\sigma)+R_k, \qquad  |R_k|\leq C'\vol(\Delta)k^{n-2},
\] 
for some $C'>0$. From this, one derives
 \begin{align*}
\left( \frac{\vol(k\Delta)-\chi(k\Delta)}{\chi(k\Delta)\vol(k\Delta)} \right) \int_{k\Delta}f(\bs x) d\mu     
=&-\left( \frac{k^{n-1}}{2}\vol(\dd \Delta)+R_k\right)\frac{1}{\chi(k\Delta)\vol(k\Delta)}\int_{k\Delta}f(\bs x){d\mu}\\
\geq&-\frac{1}{\chi(k\Delta)}\left(\frac{\vol(\partial (k\Delta))}{2\vol(k\Delta)}+C'k^{-2}\right)\int_{k\Delta}f(\bs x){d\mu}.
 \end{align*}
Next, since $(X,L)$ is uniformly K-stable, there exists $\lambda>0$ such that for any normalized convex function $f:k\Delta\rightarrow \mathbb R$, 
 \[
 \int_{\partial (k\Delta)}f(\bs x)d\sigma-\dfrac{\vol(\partial (k\Delta))}{\vol(k\Delta)}\int_{k\Delta}f(\bs x){d\mu}\geq \lambda\int_{\partial (k\Delta)}f(\bs x)d\sigma.
 \] 
Equivalently, 
\[
-\frac{\vol(\partial (k\Delta))}{\vol(k\Delta)}\int_{k\Delta}f(\bs x) d\mu\geq (\lambda-1)
\int_{\partial (k\Delta)}f(\bs x)d\sigma.
\] 
Let $C_1:=\max\set{C,C'}$. Combining the above estimate, we obtain
 \begin{align*}
\eqref{eq:Diff} \geq&  \dfrac{1}{2\chi(k\Delta)}\sum_{\bs p\in \partial(k\Delta)\cap \mathbb Z^n}f(\bs p)- \frac{M'}{k\chi(k\Delta)}\sum_{\bs p\in \partial(k\Delta)\cap \mathbb Z^n}f(\bs p)-\frac{C_1}{k^2\chi(k\Delta)}\int_{k\Delta}f(\bs x){d\mu}\\
 &-\frac{1-\lambda}{2\chi(k\Delta)}\int_{\partial (k\Delta)}f(\bs x)d\sigma.
  \end{align*}
By Corollary \ref{boundary estimates}, since each $C(k\bs p_{\alpha})$ is small, for all sufficiently large $k$, there exists $M_2>0$ independent of $k$, such that
  \begin{align*}
 \left(1+\dfrac{M_2}{k}\right) \sum_{\bs p\in \partial (k\Delta) \cap \mathbb Z^n}f(\bs p)  \geq&
 \int_{\partial (k\Delta)}f(\bs x)d\sigma. 
  \end{align*}
Putting everything together, we deduce that for all sufficiently large $k$,

\begin{align*}
\eqref{eq:Diff}  \geq&~  \dfrac{\lambda}{2\chi(k\Delta)}\sum_{\bs p\in \dd(k\Delta)\cap \bbZ^n}f(\bs p)- \dfrac{M'}{k\chi(k\Delta)}\sum_{\bs p\in \dd(k\Delta)\cap \bbZ^n}f(\bs p)-\frac{C_1}{k^2\chi(k\Delta)}\int_{k\Delta}f(\bs x){d\mu}\\
& -\frac{1-\lambda}{2\chi (k\Delta)}\left( 1+\frac{M_2}{k}\right)\sum_{\bs p\in \partial(k\Delta)\cap \mathbb Z^n}
f(\bs p)\\
=&\frac{1}{2\chi(k\Delta)}\left(\lambda -\frac{2M'+M_2(1-\lambda)}{k}\right)\sum_{\bs p\in \partial(k\Delta)\cap \mathbb Z^n}f(\bs p)
-\frac{C_1}{k^2\chi(k\Delta)}\int_{k\Delta}f(\bs x) d\mu\\
\geq&\frac{1}{2\chi (k\Delta)}
\left[
\left(\lambda -\frac{\widetilde M}{k}\right)\sum_{\bs p\in \partial(k\Delta)\cap \mathbb Z^n}f(\bs p)
+\frac{2C_1(\lambda -1)}{ka}\int_{\partial(k\Delta)}f(\bs x)d\sigma
\right],
\end{align*}
where 
\[
\widetilde M:=2M'+M_2(1-\lambda), \qquad
a=\frac{\vol(\partial \Delta)}{\vol(\Delta)}.
\] 
This concludes that for any normalized convex function $f$,
 \[
\frac{1}{\chi(k\Delta)}\sum_{\bs p\in k\Delta \cap \bbZ^n}f(\bs p)-\dfrac{1}{\vol(k\Delta)}\int_{k\Delta}f(\bs x){d\mu}\geq 0
\] 
for all sufficiently large $k$.

Since equality holds for constant functions, the same inequality holds for all convex functions $f$.
This is precisely the condition for asymptotic Chow polystability of $(X,L)$.
\end{proof}

%Sec5
\section{Semi-Canonical Triangulation}\label{sec: sct}
In this section, we describe the construction of an ideal triangulation. The first step is to triangulate $\bbR^n$.

We adapt an idea from Algebraic Topology \cite[p. $112$]{Hatcher} for this construction.
Let 
\[
I^1=[0,1]
\] 
be the unit interval in $\bbR$. 

We proceed inductively. Suppose that 
\[
S=\conv\set{\bs p_0,\dots, \bs p_k}\subset I^k:=\underbrace{[0,1]\times \dots \times [0,1]}_{\text{$k$-th}}
\] 
is a simplex in the triangulation of $I^k$. 
Consider the set 
\[
\Set{(\bs p_0,0), \dots, (\bs p_k,0), (\bs p_0,1), (\bs p_1,1), \dots, (\bs p_k,1)}.
\] 
For each $i=0,\dots, k$, define a simplex in $I^{k+1}$ by
\[
S_i^+:=\conv\Set{(\bs p_i,0),\dots, (\bs p_k,0), (\bs p_0,1),\dots, (\bs p_i,1)}.
\] 
Inductively, this defines a simplex triangulation of $I^n$ for every $n \in \mathbb N$. 

By parallel translation, this induces a type F triangulation of $\bbR^n$.
We call it the {\emph{standard triangulation}} of $\mathbb R^n$ and denote by $\mathcal T_{\mathbb R^n}$.

For example, in dimension two, this triangulation is obtained by restricting the triangulation by the lines 
\[
x=n_1, \qquad y=n_2, \qquad x+y=m, \qquad n_1, n_2, m\in \mathbb Z.
%\Set{\set{x=n_{1,0}}, \{x_2=n_{2,0}\}, \{x_1+x_2=n_{1,1}\}|n_{1,0},n_{2,0},n_{1,1}\geq 0\}.
\]
\begin{figure}[h!]
\begin{tikzpicture}[x=.5cm,y=.5cm]
\draw[black] (-1.5,-1) -- (2.5,-1);
\draw[black] (-1.5,0) -- (2.5,0);
\draw[black] (-1.5,1) -- (2.5,1);
\draw[black] (-1.5,2) -- (2.5,2);
\draw[black]   (-1,-1.5) -- (-1,2.5);
\draw[black] (0,-1.5) -- (0,2.5);
\draw[black] (1,-1.5) -- (1,2.5);
\draw[black] (2,-1.5) -- (2,2.5);
\draw[black] (1.5,-1.5) -- (-1.5,1.5);
\draw[black] (2.5,-1.5) -- (-1.5,2.5 );
\draw[black] (-1.5,0.5) -- (0.5,-1.5 );
\draw[black] (2.5,-0.5) -- (-0.5,2.5 );
\draw[black] (2.5,0.5) -- (0.5,2.5 );
\draw[black] (2.5,1.5) -- (1.5,2.5 );
\draw[black] (-1.5,-0.5) -- (-0.5,-1.5 );
\draw[fill=blue] (0,0) circle (2pt) ;
\end{tikzpicture}
\caption{Triangulation of $\bbR^2$.  
 }\label{fig: triangulation of Rn}
\end{figure}

In dimension three, the cube $I^3$ is triangulated into the following simplices. 
\begin{itemize}
\item The three $3$-simplices arising from $\conv\set{\bs e_0,\bs e_1, \bs e_2}$, where 
$\bs e_0$ denotes the origin:  
\[
\conv\set{\bs e_0, \bs e_1, \bs e_2, \bs e_3}, 
\conv\set{\bs e_1, \bs e_2, \bs e_3, \bs e_1+ \bs e_3}, 
\conv\set{\bs e_2, \bs e_3, \bs e_1+ \bs e_3, \bs e_2+ \bs e_3}.
\]
\item The three 3-simplices arising from $\conv\set{\bs e_1, \bs e_2, \bs e_1+ \bs e_2}$:  
\begin{align*}
\conv&\set{\bs e_1, \bs e_2, \bs e_1+ \bs e_2, \bs e_1+ \bs e_3}, 
\conv\set{\bs e_2, \bs e_1+\bs e_2, \bs e_1+\bs e_3, \bs e_2+\bs e_3}, \\
\conv&\set{\bs e_1+\bs e_2, \bs e_1+\bs e_3, \bs e_2+\bs e_3, \bs e_1+\bs e_2+ \bs e_3}.
\end{align*}
\end{itemize}
In the following lemma, we show that any standard triangulation is locally centrally symmetric.
The idea of the proof is the following. Suppose $S$ is an $n$-dimensional simplex in the standard triangulation $\mathcal T_{\mathbb R^n}$, and let $\bs p$ be a vertex of $S$. Then there exists an $(n-1)$-dimensional simplex 
$\underline{S}$ such that $S$ is obtained from $\underline{S}$ via the inductive construction above. By the induction hypothesis, 
the reflected simplex $\underline{S}^-(\bs p)$ with respect to $\bs p$ also belongs to $\mathcal T_{\mathbb R^{n-1}}$. Suppose that $S$ is the $i$-th simplex in $\underline S\times[0,1]$, where $0\leq i \leq n$. 
We will show that the $(n-i)$-th simplex in either $\underline{S}^-(\bs p)\times[-1,0]$ or $\underline{S}^-(\bs p)\times[1,2]$ coincides with the reflected simplex $S^-(\bs p)\in \mathcal T_{\mathbb R^n}$ (see, Example \ref{ex:Reflected}).
%Lem
\begin{lemma}\label{lem: lcs}
The standard triangulation $\mathcal T_{\bbR^n}$ is locally centrally symmetric.
\end{lemma}
\begin{proof}
We prove the claim by induction on $n$.

Since the triangulation depends on the ordering of vertices, we introduce the following notation.

\begin{enumerate}
\item For a point $\bs p \in \mathbb R^n$, define 
\[
p_i^-(\bs p):=-(\bs p_i- \bs p)+\bs p=-\bs p_i+2 \bs p.
\]
\item Let 
\[
S=\conv\set{\bs p_0,\dots, \bs p_n}\subset \bbR^n,
\]
where the vertices are ordered. For a vertex $\bs p\in \mathcal V(S)$, define 
\[
S^-(\bs p):=\conv\set{p_0^{-}(\bs p),\cdots, p_n^{-}(\bs p)}=-(S-\bs p)+\bs p.
\]
\end{enumerate} 
Our goal is to show that for every simplex $S$ in the standard triangulation and every vertex $\bs p\in \mathcal V(S)$,
the reflected simplex $S^-(\bs p)$ also belongs to the triangulation.

%Dim=1
The statement is clear in dimension one. Indeed, if $S=[\bs p, \bs p+1]$, then 
\[
S^-(p)=[\bs p-1, \bs p], \qquad S^-(\bs p+1)=[\bs p+1, \bs p+2],
\]
and both intervals belong to the standard triangulation.

Assume now that the statement holds in dimension $k$. We prove it in dimension $k+1$.
Let $S$ be a simplex in the standard triangulation, and let $\bs p\in \mathcal V(S)$ be a vertex of $S$.
Let $\bs p'=(\bs p,c)\in \mathbb R^{k+1}$, and let
\[
\widetilde{S}=\conv\set{\bs p_0',\cdots, \bs p_{k+1}'}(=:S)
\]
be a simplex containing $\bs p'$. For simplicity, we continue to denote this simplex $\widetilde{S}$ by $S$.
By construction, either the first vertex $\bs p_0'$ or the last vertex $\bs p_{k+1}'$ has last coordinate $c$.

\vskip 5pt

\noindent{\bf{Case (1).}} Suppose that the last coordinate of $\bs p_0'$ is $c$. Assume that $\bs p'$ is the $i$-th vertex, where $i\leq k$.
Let $r<k$ be the integer such that $\bs p_j'=(\bs p_j,c)$ for $0\leq j \leq r$.
By construction, the simplex $S$ is obtained from the ordered vertices
\begin{align*}
\bigl \{(\bs p_{r+1},c)& \dots, (\bs p_k, c), (\bs p_0, c),\dots, (\bs p_r,c), \\
& \qquad \qquad  (\bs p_{r+1}, c+1), \dots, (\bs p_k,c+1),  (\bs p_0,c+1),\dots, (\bs p_r,c+1) \bigr \}.
\end{align*} 
Reorder the vertices $\bs p_{r+1}, \dots, \bs p_k, \bs p_0, \dots, \bs p_r $ as $\hat {\bs p}_0, \dots, \hat {\bs p}_k$, so that
$\bs p=\hat {\bs p}_{k-r+i}$. Define
% by reordering $\bs p_{r+1}, \dots, \bs p_k, \bs p_0, \dots, \bs p_r $, where $\bs p=\hat {\bs p}_{k-r+i}$. Let
\[
\underline{S}:=\conv\Set{\hat{\bs p}_0,\cdots, \hat{\bs p}_k}.
\]
Then 
\[
S=\conv\Set{(\hat{\bs p}_{k-r},c),\dots, (\hat{\bs p_k},c),(\hat{\bs p}_0,c+1),\dots, (\hat{\bs p}_{k-r},c+1)}.
\]
By the induction hypothesis,  
\[
\underline{S}^-(\bs p):=\conv\Set{\hat{p}_k^-(\bs p),\dots, \hat{p}_{0}^-(\bs p)}\in \mathcal T_{\mathbb R^k},
\] 
where $\hat{p}_i^-(\bs p)=-(\hat{\bs p}_i-\bs p)+\bs p$, i.e., $\underline{S}^-(\bs p)$ belongs to the standard triangulation. %i.e., $\underline{S}^-(\bs p) \in \cT_{\bbR^{k}}$.

Set $\bs q_i:=\hat{p}_{k-i}^-(\bs p)$ for $0\leq i \leq k$. Applying the construction of the standard triangulation yields
\[
S'':=\conv\Set{(\bs q_{r},c-1), \dots, (\bs q_k, c-1), (\bs q_0, c), \dots, (\bs q_{r},c)}\in \cT_{\bbR^{k+1}}.
\]
Substituting the definition of $\bs q_i$, we obtain
\begin{align*}
S''=\conv\bigl\{& (\hat{p}_{k-r}^-(\bs p),c-1), (\hat{p}_{k-r-1}^-(\bs p),c-1), \dots, (\hat{p}_{0}^-(\bs p),c-1),\\
&    (\hat{p}_k^-(\bs p), c), \dots, (\hat{p}_{k-r}^-(\bs p),c) \bigr \}.
%=&S^-(\bs p),
\end{align*}
This is precisely the reflected simplex $ S^-(\bs p')$, where $\bs p'=(\hat{p}_{k-r+i}^-(\bs p),c)$. Moreover
 $\bs p'$ appears in the $k+2-(i+1)=(k+1-i)$-th position. 
Consequently, $S^-(\bs p')\in\cT_{\bbR^{k+1}}$.

\vskip 5pt

\noindent{\bf{Case (2).}} 
Suppose instead that the last coordinate of $\bs p_{k+1}'$ is $c$. Assume that $S$ is generated from a suitable subset of 
\begin{align*}
\bigl\{ &(\bs p_{r+1},c-1) \dots, (\bs p_k, c-1), (\bs p_0, c-1),\dots, (\bs p_r,c-1), \\
 & (\bs p_{r+1}, c), \dots, (\bs p_k,c), (\bs p_0,c), \dots, (\bs p_r,c) \bigr \}.
\end{align*}
%Similar to the previous case, we let 
Define the simplex
\[
\underline{S}:=\conv\Set{\bs p_{r+1}, \dots, \bs p_k, \bs p_0,\dots, \bs p_r}=\conv\Set{\hat{\bs p}_0,\dots, \hat{\bs p}_k},
\] 
where $\bs p=\hat{\bs p}_{i-r-1}$. Then %(The first point is at the $0$-th position). Also, 
\[
S=\conv\Set{(\hat{\bs p}_{k-r},c-1),\dots, (\hat{\bs p}_k,c-1),(\hat{\bs p}_0,c),\dots, (\hat{\bs p}_{k-r},c)}.
\]
Again, by the induction hypothesis, %we see that 
\[
\underline{S}^-(\bs p):=\conv\Set{\hat{ p}_k^-(\bs p),\dots, \hat{ p}_0^-(\bs p)}\in \cT_{\bbR^{k}}.
\] 
The construction then yields 
\begin{align*}
S^-(\bs p')=\conv\bigl \{& (\hat{p}_{k-r}^-(\bs p),c), (\hat{p}_{k-r-1}^-(\bs p),c), \dots, (\hat{p}_{0}^-(\bs p),c), \\
& (\hat{p}_k^-(\bs p), c+1), \cdots, (\hat{p}_{k-r}^-(\bs p),c+1)\bigr \}, 
\end{align*}
where $\bs p'=(\hat{p}^-_{i-r-1}(\bs p),c)$. Therefore, $S^-(\bs p') \in \mathcal T_{\mathbb R^{k+1}}$.

Combining two cases completes the induction. Hence the standard triangulation is locally centrally symmetric.
\end{proof}

%Ex
\begin{example}\label{ex:Reflected}\rm 
Let $\bs p=(0,0)$ and $\bs p'=(0,0,0)$. Suppose that
\[
S=\conv\Set{(-1,0,0), (0,0,0),(-1,1,0), (-1,0,1)}.
\] 
Projecting onto the first two coordinates, we obtain
\[
\underline{S}=\conv\Set{(-1,0), (0,0),(-1,1)}.
\] 
Then the reflected simplex with respect to $\bs p$ is
\[
\underline{S}^-(\bs p)=\conv\Set{(1,-1), (0,0),(1,0)}.
\] 
Therefore, lifting back to dimension three, we obtain
\[
S^-(\bs p')=\conv\Set{(1,0,-1),(1,-1,0), (0,0,0),(1,0,0)}.
\]
\begin{figure}[h!]
\begin{tikzpicture}[x=.5cm,y=.5cm]
\draw[black] (-1,0) -- (0,0);
\draw[black] (0,0) -- (-1,1);
\draw[black] (-1,1) -- (-1,0);
\draw[red] (1,0) -- (0,0);
\draw[red]   (0,0) -- (1,-1);
\draw[red] (1,-1) -- (1,0);

\draw[fill=blue] (0,0) circle (2pt) ;
\end{tikzpicture}
\caption{$\underline{S}$ and $\underline{S}^-(p)$ }
\end{figure}

\begin{figure}[h!]
\tdplotsetmaincoords{20}{10}
\begin{tikzpicture}[tdplot_main_coords]

\draw[black, thick]  (-1,0,0) -- (0,0,0) -- (-1,1,0) -- (-1,0,0) ;
\draw[black, thick] (-1,0,0) -- (0,0,0)-- (-1,0,1)-- (-1,0,0) ;
\draw[black, thick] (0,0,0) -- (-1,1,0)-- (-1,0,1)-- (0,0,0) ;

\draw[red, thick]  (1,0,0) -- (0,0,0) -- (1,-1,0) -- (1,0,0) ;
\draw[red, thick] (1,0,0) -- (0,0,0)-- (1,0,-1)-- (1,0,0) ;
\draw[red, thick] (0,0,0) -- (1,-1,0)-- (1,0,-1)-- (0,0,0) ;

\end{tikzpicture}
\caption{$S$ and $S''$ }
\end{figure}
\end{example}
%Ex2
\begin{example}\rm
We illustrate the proof of Lemma \ref{lem: lcs} with an explicit example. %for $k=2$, $r=1$, $c=1$.
Let 
\[
S=\conv\Set{(1,0,0),(0,1,0),(0,0,1),(1,0,1)},
\]
and let $\bs p'=(1,0,1)$. We construct the reflected simplex $S^-(\bs p')$.
%In this case, we see that $k=2$, $r=1$, $c=1$.

Observe that $S$ is the convex hull of the second, third, fourth, and fifth points of the set
\[
\Set{(0,0,0), (1,0,0), (0,1,0),(0,0,1),(1,0,1), (0,1,1)}.
\]
Projecting onto the first two coordinates gives
\[
\underline{S}=\Set{(1,0), (0,1),(0,0)}.
\]
Reordering the vertices, we write $\hat{\bs p}_0=(0,0), \hat{\bs p}_1=(1,0), \hat{\bs p}_2=(0,1)$.
In this case, $k=2$, $r=1$, $c=1$, and
\[
S=\Set{( \hat{\bs p}_1,0),( \hat{\bs p}_2,0),( \hat{\bs p}_0,1),( \hat{\bs p}_1,1)}.
\]
Next, reflect $\underline{S}$ with respect to $\bs p=(1,0)$. We obtain
\[
\underline{S}^-(\bs p)=\conv\Set{(2,-1),(1,0),(2,0)}=\conv\Set{\hat{ p}_2^-(\bs p),\hat{ p}_1^-(\bs p), \hat{ p}_0^-(\bs p)}.
\]
Applying the construction of the standard triangulation yields
\begin{align*}
S^-(\bs p')&=\conv\Set{(\hat{ p}_1^-(\bs p),1),(\hat{ p}_0^-(\bs p),1),(\hat{ p}_2^-(\bs p),2),(\hat{ p}_1^-(\bs p),2)}\\
&=\conv\Set{(1,0,1),(2,0,1),(2,-1,2),(1,0,2)}.
\end{align*}
On the other hand,
\[
S-\bs p'=\conv\Set{(0,0,-1),(-1,1,-1),(-1,0,0),(0,0,0)},
\]
and hence
\[
-(S-\bs p')+\bs p'=\conv\Set{(1,0,1),(2,-1,2),(2,0,1),(1,0,2)}=S^-(\bs p').
\]
Thus, the construction indeed gives the reflected simplex.
\end{example}
Recall that for a triangulation $\mathcal T$ and a lattice point $\bs p\in \mathbb Z^n$, we denote by $n(\bs p)$  the number of simplices in $\mathcal T$ containing $\bs p$ as a vertex.
%Lem
\begin{lemma}\label{lem: n(p) of standard triangulation}
 For any $\bs p\in \bbZ^n$, the standard triangulation $\mathcal T_{\mathbb R^n}$ is a simplex triangulation of $\bbR^n$ satisfying $n(\bs p)=(n+1)!$.

\end{lemma}
\begin{proof}
From the definition $\mathcal T_{\mathbb R^n}$ is a simplex triangulation. Hence it suffices to show that $n(\bs p)=(n+1)!$.
We prove the statement by induction on $n$.

For $n=1$, this is obvious. 
For $n=2$, let $T(I^2)$ be a triangulation of $I^2$. Then $T(I^2)$ 
consists of two triangles. By taking parallel transformations of $T(I^2)$ around a lattice point $\bs p$, we obtain the standard triangulation $\mathcal T_{\bbR^2}$, and each lattice point is contained in $3!=6$ triangles. 

Assume now that the statement holds in dimension $k$. We prove it in dimension $k+1$. Let
\[
S_k=\conv\set{\bs p_0,\dots, \bs p_k}
\]
be a simplex in $\mathcal T_{\mathbb R^k}$. %the standard triangulation of $\mathbb R^n$. 
The prism $S_k\times [0,1]$ is triangulated into simplices
% observe that for any simplex $S_n=\conv\set{\bs p_0,\dots, \bs p_n}$, when we consider $S_n\times [0,1]$, this is triangulated by 
\[
S_{k+1}^i=\conv\Set{(\bs p_i,0),\dots, (\bs p_k,0),(\bs p_0,1),\dots, (\bs p_i,1)},
\] 
for $0\leq i\leq k$. 
Suppose that $\bs p=\bs p_i\in S_k$.
Then there are exactly $i+1$ simplices in $S_k\times [0,1]$ containing 
$(\bs p,0)$. 

Similarly, there are exactly $k-i+1$ simplices in  $S_k\times [-1,0]$  containing $(\bs p,0)$. Therefore, by the induction hypothesis, 
\[
n(\bs p)=[(i+1)+(k-i+1)](k+1)!=(k+2)!.
\] 
Hence the statement holds in dimension $k+1$.
\end{proof}

%Example
\begin{example}\rm
%For example, 
When $n=1$, the induction step from dimension one to dimension two can be visualized as follows.
% goes from 1 dimension to 2 dimensions. 

 \begin{figure}[h!]
\begin{tikzpicture}[x=.5cm,y=.5cm]
\draw[thick] (1,0) -- (0,1);
\draw[thick] (0,1) -- (-1,1);
\draw[thick] (0,1) -- (-1,1);
\draw[thick] (-1,1) -- (-1,0);
\draw[thick] (-1,0) -- (0,-1);
\draw[thick] (0,-1) -- (1,-1);
\draw[thick] (1,-1) -- (1,0);
\draw[red,thick] (1,0) -- (-1,0);
\draw[red,thick] (0,1) -- (0,-1);
\draw[red,thick] (-1,1) -- (1,-1);
\end{tikzpicture}
\caption{The triangulation near (0,0)}
\end{figure}

There are exactly two intervals containing $0$, namely $[0,1]$ and $[-1,0]$. 
For $S=[0,1]$, the prism $S\times [0,1]$ contains one triangle containing $(0,0)$ (corresponding to the case $i=0$), while $S\times [-1,0]$ contains two such triangles.

Similarly, for $S=[-1,0]$, there are two triangles containing $(0,0)$ in $S\times [0,1]$, and one triangle containing
$(0,0)$ in $S\times [-1,0]$.

Altogether, we obtain
\[
1+2+2+1=6=3!
\]
triangles containing $(0,0)$.
\end{example}

%Claim
\begin{claim}\label{claim: all planes}
Let $H$ be a hyperplane arising as the affine span of a codimension-one face of a simplex in $\mathcal T_{\mathbb R^n}$, passing through the origin. 
Then $H$ is one of the hyperplanes
\begin{align*}
\Set{x_i+x_{i+1}+\dots +x_j=0| 1\leq i\leq j\leq n}.
\end{align*}
\end{claim}
\begin{proof} 
We first show that the triangulation $\mathcal T_{\mathbb R^n}$ contains all hyperplanes of the above form. We proceed by induction on the dimension. 

For $n=2$, the hyperplanes passing through the origin are 
\[
\Set{x_1=0},\qquad  \Set{x_2=0}, \qquad \set{x_1+x_2=0}.
\]

\begin{figure}[h!]
\begin{tikzpicture}[x=.5cm,y=.5cm]
\draw[thick] (1,0) -- (0,1);
\draw[thick] (0,1) -- (-1,1);
\draw[thick] (0,1) -- (-1,1);
\draw[thick] (-1,1) -- (-1,0);
\draw[thick] (-1,0) -- (0,-1);
\draw[thick] (0,-1) -- (1,-1);
\draw[thick] (1,-1) -- (1,0);
\draw[red,thick] (1,0) -- (-1,0);
\draw[red,thick] (0,1) -- (0,-1);
\draw[red,thick] (-1,1) -- (1,-1);
\end{tikzpicture}
\caption{The triangulation at the center of $(0,0)$}
\end{figure}

Assume that the statement holds in dimension $n$. Observe that if a face $F$ of the boundary $ \partial S$ of a simplex $S\in \mathcal T_{\mathbb R^n}$ is contained in a hyperplane
\[
\Set{\bs x=(x_1, \dots, x_n)\in \mathbb R^n| f(\bs x)=c},
\]
then the corresponding face of the prism $S\times [0,1]\subset \mathbb R^{n+1}$, % is
\[
{F}'=\Set{(\bs x,0)\in \mathbb R^{n+1}| \bs x\in F},
\]
is contained in the same hyperplane after identifying $\mathbb R^n\cong\set{x_{n+1}=0}\subset \mathbb R^{n+1}$.

Since the number of simplices grows rapidly in higher dimensions, the inductive step is not intended to determine all supporting hyperplanes explicitly, but only to show that all hyperplanes listed in the statement occur in the triangulation.
Thus it suffices to show that the triangulation contains hyperplanes of the form
\[
x_i+x_{i+1}+\cdots+x_{n+1}=c_i
\]
for suitable constants $c_i$. By parallel translation, this then yields hyperplanes  
\[
x_i+x_{i+1}+\cdots+x_{n+1}=0, \qquad i=1,\dots, n+1.
\] 
Consider the prism
\[
S_n\times [0,1], \qquad S_n=\conv\set{\bs 0,\bs e_1,\dots, \bs e_n}\subset \mathbb R^n.
\]
By construction, $S_n\times [0,1]$ is triangulated into $n+1$ simplices obtained as the convex hulls of $n+2$ consecutive points in
\[
\Set{\bs 0, \bs e_1,\dots, \bs e_n, \bs e_{n+1}, \bs e_1+\bs e_{n+1}, \dots, \bs e_{n}+\bs e_{n+1}}\subset \mathbb R^{n+1}.
\] 
Writing $\bs e_0=\bs 0\in \mathbb R^{n+1}$, we see that these simplices are precisely
\[
\Delta_i:=\conv\Set{\bs e_i,\dots, \bs e_n, \bs e_0+\bs e_{n+1},\bs e_1+\bs e_{n+1}, \dots, \bs e_i+\bs e_{n+1}},
\quad i=0,\dots, n.
\]

We now analyze the boundary faces of $\Delta_i$. We claim that for each $i=1,\cdots, n$,
one of the faces of $\Delta_i$ is contained in a hyperplane
\[
x_i+\cdots+x_{n+1}=c_i,
\] 
while one of the faces of $\Delta_0$ is contained in $x_{n+1}=0$.
Indeed, for
\[
\Delta_0=\conv\set{\bs e_0, \bs e_1,\dots, \bs e_{n+1}},
\]  
the face spanned by the first $n+1$ vertices lies in the hyperplane $x_{n+1}=0$.

For $\Delta_i$ $(i\geq 1)$, the relevant face of $\Delta_i$ is spanned by 
\[
\Set{\bs e_i,\dots, \bs e_n, \bs e_0+\bs e_{n+1},\bs e_1+\bs e_{n+1}, \dots, \bs e_{i-1}+\bs e_{n+1}}.
\]  
To determine the supporting hyperplane of this face, subtract the vertex $\bs e_i$ to obtain tangent vectors
\[
\Set{\bs e_{i+1}-\bs e_i,\dots, \bs e_n-\bs e_i, \bs e_{n+1}-\bs e_i,\bs e_1+\bs e_{n+1}, \dots, \bs e_{i-1}+\bs e_{n+1}-\bs e_i}.
\] 
Then for all $r\geq 1$,
\begin{align*}
\langle \bs e_i+\dots+\bs e_{n+1}, \bs e_{i+r}-\bs e_{i}\rangle&=|\bs e_{i+r}|^2-|\bs e_i|^2=0, \quad \text{and}\\
\langle \bs e_i+\dots+\bs e_{n+1}, \bs e_{i-r}+\bs e_{n+1}-\bs e_{i}\rangle&=|\bs e_{n+1}|^2-|\bs e_i|^2=0.
\end{align*}
Hence, $\bs e_i+\dots+\bs e_{n+1}$ is a normal vector to the hyperplane containing this face. Therefore the hyperplane is given by
\[
x_i+\dots+x_{n+1}=c_i.
\]
After a parallel translation, we obtain a simplex whose face lies in
\[
x_i+\cdots+x_{n+1}=0.
\]
This proves that all hyperplanes listed in the statement occur as supporting hyperplanes of simplices meeting the origin.

%(2nd Half)
Next, we show that these hyperplanes divide a neighborhood of the origin into exactly $(n+1)!$ chambers. By Lemma \ref{lem: n(p) of standard triangulation}, this implies that there are no additional supporting hyperplanes through the origin.

Define
\[
q_0:=0, \qquad q_i:=x_1+\dots + x_i.
\]
Then every hyperplane in Claim $\ref{claim: all planes}$ is of the form
\[
q_j-q_i=0 \qquad (0\leq i<j \leq n).
\]
Hence the complement of these hyperplanes is precisely the set where
\[
q_{\sigma(0)}>q_{\sigma(1)}>\dots >q_{\sigma(n)}
\]
for some permutation $\sigma \in \mathfrak S_{n+1}$. Therefore the arrangement divides a neighborhood of the origin into exactly
$(n+1)!$ chambers.

By Lemma \ref{lem: n(p) of standard triangulation}, exactly $(n+1)!$ simplices of the standard triangulation meet at the origin, and each simplex determines a unique chamber bounded by supporting hyperplanes through the origin.
Since the above hyperplanes already produce exactly $(n+1)!$ chambers, there can be no additional supporting hyperplanes through the origin.

Finally, since the standard triangulation is locally centrally symmetric and globally obtained by parallel translations from the triangulation of $I^n=[0,1]^n$,
there are no aditional hyperplanes anywhere in the standard triangulation. 
Therefore the supporting hyperplanes through the origin are precisely the hyperplanes 
\[
\Set{x_i+x_{i+1}+\dots +x_j=0| 1\leq i\leq j\leq n}.
\]
\end{proof}

Let $\Pi_n$ denote the standard cone:
\[
\Pi_n:=\Set{\bs x=(x_1,\cdots, x_n)\in \bbR^n | x_i\geq 0 \text{ for } 1\leq i\leq n}.
\]

%Lem
\begin{lemma}\label{lem:Lee25}
Under the standard triangulation $\mathcal T_{\mathbb R^n}$, %of $\bbR^n$, 
we have
$n(\bs p)=(n+1)!$ for all $\bs p\in \bbZ^n$. 
Moreover, this triangulation induces on $\Pi_n$ the restriction of the locally centrally symmetric triangulation of $\mathbb R^n$.
Furthermore,
\[
n(\bs p)=\dfrac{(n+1)!}{(k+1)!} \quad \text{ for all } \quad \bs p\in \mathrm{Skel}_{n-k}(\Pi_n)^o\cap \bbZ^n.
\]  
\end{lemma}
\begin{proof}
By Lemma \ref{lem: lcs} and Lemma \ref{lem: n(p) of standard triangulation}, it remains only to show that
\[
n(\bs p)=\frac{(n+1)!}{(k+1)!} 
\]  
for every lattice point $\bs p\in \mathrm{Skel}_{n-k}(\Pi_n)^o\cap \bbZ^n$.

Take $\bs p=(p_1, \dots, p_n)\in \mathrm{Skel}_{n-k}(\Pi_n)^o\cap \bbZ^n$. Without loss of generality, we may assume that
the triangulation is centered at the origin. We then apply Claim \ref{claim: all planes} to enumerate all supporting hyperplanes passing through the
origin.

Define
\[
q_0=0, \qquad  q_i:=p_1+\dots +p_i.
\]
By Claim \ref{claim: all planes}, the supporting hyperplanes through the origin are precisely
\[
q_j-q_i=0 \qquad (0\leq i<j\leq n).
\]
Since $\bs p \in \mathrm{Skel}_{n-k}(\Pi_n)^o$, exactly $k$ coordinates of $\bs p$ vanish.
Equivalently, among the numbers $q_0, q_1, \dots, q_n$, there are exactly $k+1$ equal consecutive values.

Therefore, the number of chambers adjacent to $\bs p$ equals the number of total orderings of $q_0, \dots, q_n$ modulo permutations of these $k+1$ equal values: $\frac{(n+1)!}{(k+1)!}$.
Hence, 
\[
n(\bs p)=\frac{(n+1)!}{(k+1)!}
\]
for all $\bs p\in \mathrm{Skel}_{n-k}(\Pi_n)^o\cap \bbZ^n$. This completes the proof.
\end{proof}

%Remark
\begin{remark}\rm
The same argument shows that for the induced triangulation on the dilated simplex $ mS_n$ with $m\in \mathbb N$,    
\[
n(\bs p)=\dfrac{(n+1)!}{(k+1)!} 
\]  
for all $\bs p \in \mathrm{Skel}_{n-k}(mS_n)^o\cap \mathbb Z^n$.
\end{remark}
%Lem

%Cor
\begin{corollary}\label{standard cone is small}
The standard cone $\Pi_n$ is small with respect to the standard triangulation on $\Pi_n$ and the triangulation on $\dd\Pi_n$ induced from the standard triangulation on each face.       
\end{corollary}
\begin{proof}
By Lemma \ref{lem:Lee25}, for every point $\bs p\in \Pi_n^o\cap\bbZ^n$, we have
\[
n(\bs p)=(n+1)!,
\] 
while for every point $\bs p\in \dd\Pi_n\cap\bbZ^n$, we have
\[
n(\bs p)\leq\dfrac{(n+1)!}{2}.
\] 
Thus it remains to show that for every $\bs p\in \dd\Pi_n\cap\bbZ^n$, we have
\[
m(\bs p)\leq n!.
\] 
Suppose that $\bs p\in \mathrm{Skel}_{n-k}(\Pi_n)^o\cap \bbZ^n$. Then $\bs p$ lies in the intersection of exactly $k$ facets, i.e.,
 \[
 \bs p\in 
 \smallcap_{\alpha=1}^k(\dd\Pi_n)_{i_\alpha}.
 \] 
Moreover, for each facets $(\partial \Pi_n)_{i_\alpha}$, 
\[ 
\bs p\in\mathrm{Skel}_{n-1-(k-1)}\bigl( (\dd \Pi_n)_{i_\alpha}\bigr)^o\cap \bbZ^n .
\] 
Each $k$ facet containing $\bs p$ is isomorphic to $\Pi_{n-1}$. Hence, by Lemma \ref{lem:Lee25}, the number of simplices in the induced boundary triangulation containing $\bs p$ inside a fixed facet equals $\frac{n!}{k!}$.
Since $\bs p$ lies on exactly $k$ facets, we obtain
\[
m(\bs p)=\dfrac{n!}{k!}\cdot k\leq n!
\] 
for all $k=1,\dots, n$. This proves the claim.
\end{proof}

The next step is to systematically construct an {\emph{ideal}} triangulation of a given vertex cone together with a type F triangulation on its boundary. 

Let 
\[
C=C(\bs 0; \bs v_1, \dots ,\bs v_r)\subset \bbR^n
\] 
be a cone with the vertex $\bs 0$ and generators $\bs v_1, \dots , \bs v_r$.
%Strategy
\begin{enumerate}
\item[Step 1.] Define $\Delta=\conv\set{\bs 0, \bs v_1, \dots, \bs v_r}$. Choose a triangulation of $\Delta$ such that every simplex $S$ in the triangulation satisfies $\vol(S)=1/n!$.
Denote this triangulation by $\mathcal T(\Delta)=\Set{S_j}_{j=1}^\alpha$.
Taking all simplices $S_j\in \mathcal T(\Delta)$ satisfying $\bs 0\in \mathcal{V}(S_j)$, define %$R$ as the union of such simplices:  
\[
R:=\smallcup_{1\leq j\leq \alpha, \atop  \bs 0\in \mathcal{V}(S_{j})}S_j.
\]
\item[Step 2.]
For each simplex $ S\in \mathcal T(\Delta)$, define the cone 
$C_S:=\Set{t\bs x | t\geq 0, \bs x\in S}$. Each cone $C_S$ is isomorphic to the standard cone. Moreover,
\[
\smallcup_{j=1}^\alpha C_{S_j}=C.
\] 

\item[Step 3.]  Apply the standard triangulation from Lemma \ref{lem:Lee25} to each cone $C_S$.   
In this way, we obtain a triangulation of $C$. Since the triangulation on each $C_S$ is strongly ideal, the resulting triangulation 
of $C$ is ideal.  

\item[Step 4.] For each facet $\dd C_{r}$ of $C$,  
consider those cones $C_{S}$ satisfies 
\[
\dim (\dd C_r\cap C_{S})=n-1.
\] 
Observe that $\dd C_r\cap C_{S}$ is isomorphic to the $(n-1)$-dimensional standard cone $\Pi_{n-1}$.
Therefore, we apply the standard triangulation to each intersection $\dd C_r\cap C_{S}$. This induces a type F triangulation on every facet $\dd C_r$, and hence on the entire boundary $\partial C$. 
\end{enumerate}
%Def
\begin{definition}\label{def:SCT}\rm
We call the resulting triangulation on $C$ and $\dd C$ the {\emph{semi-canonical triangulation}}. 
\end{definition}
Note that semi-canonical triangulations depend on the choice of the simplex triangulation of $\Delta$, and hence are not unique. 
Nevertheless, the construction uses only minimal combinatorial data and produces triangulations that are as regular as possible.

%Example
\begin{example}\rm
Consider the cone $C((0,0); (2,1), (2,-1))$. Then 
\[
\Delta=\conv\set{(0,0), (2,1),(2,-1)},
\]
shown as the green triangle in Figure \ref{fig:Q}.
%, which shown as green in the following picture. 

\begin{figure}[h!]
\begin{tikzpicture}[x=.5cm,y=.5cm]
\draw[thick] (2,1) -- (4,2);
\draw[thick] (2,-1) -- (4,-2);
\draw[red,thick] (0,0) -- (2,0);
\draw[red,thick] (1,0) -- (2,1);
\draw[red,thick] (1,0) -- (2,-1);
\draw[green,thick] (2,1) -- (2,-1);
\draw[green,thick] (2,1) -- (0,0);
\draw[green,thick] (2,-1) -- (0,0);
\end{tikzpicture}
\caption{$\Delta$ with the simplex triangulation}\label{fig:Q}
\end{figure}
As illustrated in Figure \ref{fig:R}, the %union of simplices 
region $R$ in Step 1 consists of the two triangles $S_1=\conv\set{(0,0),(1,0),(2,1)}$ and $S_2=\conv\set{(0,0),(1,0),(2,-1)}$.

\begin{figure}[h!]
\begin{tikzpicture}[x=.5cm,y=.5cm]
\draw[red,thick] (0,0) -- (1,0);
\draw[red,thick] (1,0) -- (2,1);
\draw[red,thick] (1,0) -- (2,-1);
\draw[red,thick] (2,1) -- (0,0);
\draw[red,thick] (2,-1) -- (0,0);
\end{tikzpicture}
\caption{$R$ with the simplex triangulation}\label{fig:R}
\end{figure}
Extending $R$ generates two standard cones, and we apply the standard triangulation to each cone. 
The corresponding triangulations of the two cones are shown in green and blue, respectively.
%We show the triangulation on the two cones in green and blue respectively.  
\begin{figure}[h!]
\begin{tikzpicture}[x=.5cm,y=.5cm]
\draw[thick] (6,3) -- (0,0);
\draw[thick] (6,-3) -- (0,0);
\draw[thick] (3,0) -- (0,0);
%\draw[thick] (6,3) -- (3,0);
%\draw[thick] (6,-3) -- (3,0);
\draw[red,thick] (0,0) -- (3,0);
\draw[green,thick] (1,0) -- (2,1);
\draw[blue,thick] (1,0) -- (2,-1);
\draw[green,thick] (2,1) -- (0,0);
\draw[blue,thick] (2,-1) -- (0,0);
\draw[green,thick] (2,1) -- (4,1);
\draw[blue,thick] (2,-1) -- (4,-1);
\draw[green,thick] (4,2) -- (5,2);
\draw[blue,thick] (4,-2) -- (5,-2);

\draw[green,thick] (2,0) -- (4,2);
\draw[blue,thick] (2,0) -- (4,-2);

\draw[green,thick] (1,0) -- (5,2);
\draw[blue,thick] (1,0) -- (5,-2);

%\draw[green,thick] (2,0) -- (4,1);
\draw[blue,thick] (2,0) -- (4,-1);
\end{tikzpicture}
\caption{An ideal simplex triangulation}\label{fig:3R}
\end{figure}
\end{example}
We may now reformulate Theorem \ref{theo: main theorem 1} as follows. %to the following:
%Thm
\begin{theorem}\label{theo: main theorem 2}
Let $(X,L)$ be a polarized toric variety with integral polytope $\Delta$. Suppose that: 
\begin{enumerate}
    \item $(X,L)$ is uniformly K-stable;
    \item all Futaki-Ono invariants vanish;
    \item every vertex cone of $\Delta$ admits a semi-canonical triangulation which is small.
    %under some semi-canonical triangulations are small;
\end{enumerate}
Then $(X,L)$ is asymptotically Chow polystable.
\end{theorem}
%Cor
\begin{corollary}\label{coro: K stable manifold is Chow stable}
Let $(X,L)$ be an $n$-dimensional uniformly K-stable polarized toric manifold such that all Futaki-Ono invariants vanish. 
Then $(X,L)$ is asymptotically Chow polystable.
\end{corollary}
\begin{proof}
Since $(X,L)$ is smooth, every vertex cone is isomorphic to the standard cone $\Pi_n$. By Corollary \ref{standard cone is small},
the standard cone with respect to the standard triangulation is small. Therefore, Theorem \ref{theo: main theorem 2}
implies that $(X,L)$ is asymptotically Chow polystable.
\end{proof}

%Sec6
\section{Symmetric weakly reflexive polytopes}\label{sec:SymWeak}
An integral polytope $\Delta\subset \mathbb R^n$ is called {\emph{weakly reflexive}} if there exists some constant $c$ such that 
\[
\Delta=\smallcap_{i=1}^{\alpha}\set{\ell_i(\bs x)\leq c},
\] 
where each $\ell_i(\bs x)=\sum_{j=1}^na_{i,j}x_j$ is a primitive integral linear form, i.e.,
\[
\gcd(c,a_{i,1}, \dots, a_{i,n})=1.
\]
In particular, when $c=1$, the polytope $\Delta$ is a reflexive polytope.
A weakly reflexive polytope $\Delta$ is called {\emph{symmetric}} if the action of its symmetry group
$G\subset \SL(n,\mathbb Z)$ on $\Delta$ has a unique fixed point.
Necessarily, this fixed point is the origin.
%Lem
\begin{lemma}\label{lemm lambda stability of  weakly reflexive symmetric toric variety}
 Let $\Delta$ be an $n$-dimensional symmetric weakly reflexive polytope. Then, $\Delta$ is uniformly K-stable.   
\end{lemma}
\begin{proof}
We adapt the argument of \cite[Proposition $4.6$]{Yao22}.
First, define a map $\varphi: \Delta\backslash\Set{(0, \dots ,0)} \rightarrow \dd \Delta \times (0,c]$ as follows.
For any point $\bs p\in \Delta\setminus \set{\bs 0}$, let $\bs q\in \partial \Delta$ be the unique intersection point of the ray
$\Set{t\bs p \mid t\ge 0}$ with the boundary $\partial \Delta$.
Define
\[
t=\frac{c|\bs p|}{|\bs q|}, \qquad  \varphi(\bs p):=(\bs q,t).
\] 
Let $G\subset SL(n,\bbZ)$ be the symmetry group, and let $f:\Delta \rightarrow \bbR$ be a $G$-invariant convex function. 
Since the origin is the unique fixed point of $G$, convexity implies that $f$ attains its minimum at $\bs 0$.
%Then, $f(0)$ attains the minimum. 

We now define a function $g_f: \Delta \to \bbR$ by the following: 
\[
g_f\circ \varphi^{-1}(\bs q,t):=\dfrac{t}{c}f\circ \varphi^{-1}(\bs q,1), %\quad \text{ and } \quad 
\qquad g_f(\bs 0)=0.
\]
The convexity of $f$ implies that $g_f(\bs x) \geq f(\bs x)$, whereas we see that 
\[
\int_{\dd \Delta}g_f(\bs q)d\sigma=\int_{\dd \Delta}f(\bs q)d\sigma.
\]
Now, we write the point in $\Delta$ by $(\bs x,t)\in \dd\Delta\times (0,c]$. 
Then the straightforward computation shows that
\begin{align*}
    \int_{\Delta}f(\bs x,t)d\mu\leq &\int_{\Delta}g_f(\bs x)d\mu
    %=\int_{0}^c\int_{t(\dd \Delta)}\dfrac{t}{c}\cdot g_f(\bs x,1)d\mu\\
    =&\int_0^c\dfrac{t}{c}\cdot t^{n-1}dt\int_{\dd \Delta}f(\bs x,1)d\sigma
    =\dfrac{c^n}{n+1}\int_{\dd \Delta}f(\bs x)d\sigma.
\end{align*}
Also, we have
\begin{align*}
    \int_{\Delta}d\mu=&\int_0^ct^{n-1}dt\int_{\dd \Delta}d\sigma  =\dfrac{c^n}{n}\int_{\dd \Delta}d\sigma.
\end{align*} 
Therefore, we conclude $a:=\frac{\vol(\partial \Delta, d\sigma)}{\vol(\Delta)}=\frac{n}{c^n}$. This implies that
\begin{align}
\begin{split}\label{ineq:DFInv}
    L_a(f)=&\int_{\dd \Delta}f(x)d\sigma-\dfrac{n}{c^n}\int_{\Delta}f(x)d\mu\\
    \geq& \int_{\dd \Delta}f(x)d\sigma-\dfrac{n}{c^n}\left(\dfrac{c^n}{n+1}\int_{\dd \Delta}f(x)d\sigma\right)\\
    =&\dfrac{1}{n+1}\int_{\dd \Delta}f(x)d\sigma.
\end{split}
\end{align}
Hence, $\Delta$ is uniformly K-stable.
\end{proof}

Notice that the Futaki-Ono invariant vanishes for any symmetric weakly reflexive polytope.
Therefore we obtain the following corollary. 
%Cor
\begin{corollary}\label{cor:Fano}
Let $\Delta$ be a symmetric weakly reflexive polytope. 
Suppose that every vertex cone of $\Delta$ admits a semi-canonical triangulation which is small.
 Then the associated polarized toric variety is asymptotically Chow polystable.
\end{corollary}
\begin{proof}
By Lemma \ref{lemm lambda stability of  weakly reflexive symmetric toric variety}, $\Delta$ is uniformly K-stable.
Moreover, all Futaki-Ono invariants vanish because $\Delta$ is symmetric. 
Hence, Theorem \ref{theo: main theorem 2} implies that the associated polarized toric variety is asymptotically  Chow polystable.
 \end{proof}

%Sec7
\section{Examples}\label{sec:Examples}

%Ex
\begin{example}[$D_n$, see \cite{Lee25}]\rm
Consider the polytope $\Delta_n':=\conv\Set{\pm \bs e_i}_{i=1}^n$, where $\bs e_1=(1,0, \dots ,0), \dots , \bs e_n=(0, \dots,0,1)$ denote the standard basis of $\mathbb R^n$.
Let $D_n$ be the associated projective toric variety.

First, $D_n$ is a symmetric toric Fano variety. 
Moreover, all vertex cones $C(\bs p_i)$ are isomorphic to the cone generated by the origin %$(0, \dots ,0)$ 
and the vectors 
\[
(\pm 1,0, \dots ,0,1), (0, \pm 1,0, \dots ,0,1), \dots , (0, \dots ,0,\pm 1,1).
\]
We denote this cone by $C$.  
Define
\[
\Delta=\conv\Set{(0, \dots ,0), (\pm 1,0, \dots ,0,1), \cdots , (0, \dots ,0,\pm 1,1)}.
\] 
We triangulate $\Delta$ using hyperplanes 
\[
\Set{\bs x=(x_1, \dots, x_n)\in \mathbb R^n | x_i=0}, \qquad i=1, \dots, n-1.
\]
This induces a semi-canonical triangulation on the cone $C$. We then obtain
%by $\{x_i=0\}$ for $i=1,\cdots, n-1$. We can construct a semi canonical triangulation on the cone $C$. We see that for every points,
\[
n(\bs p)=
\begin{cases}
2^{n-1} & \bs p=(0,\dots,0) \\
\dfrac{2^i(n+1)!}{(i+1)!}  & \bs p\neq \bs 0 \text{ and } \bs p\in \displaystyle  \smallcap_{r=1}^i \Set{x_{k_r}=0}, 
\\
\dfrac{(n+1)! }{2}   & \bs p\neq \bs 0 \text{ and } \bs p\in  \partial C, \\
(n+1)!  & \text{otherwise},
\end{cases}
\]
where $1\le k_1<\cdots < k_i\le n-1$.
Furthermore, if $\bs p$ is an interior point of a boundary facet,
then $\bs p$ lies in an $(n-1)$-simplex. Hence $m(\bs p)=n!$.
More generally, we obtain:
\begin{itemize}
\item if $\bs p$ lies on an $(n-2)$-dimensional face, then
\[
m(\bs p)=\frac{n!}{2!}\cdot 2=n!;
\]
\item if $\bs p$ lies on a codimension-$k$ face $(k\geq 3)$, then
\[
m(\bs p)=\dfrac{2^{k-1} n!}{k!}.
\]
\end{itemize}
Consequently, $C$ is a small cone under this semi-canonical triangulation. 
Therefore, by Corollary \ref{cor:Fano}, the varieties $D_n$ are  asymptotically  Chow polystable.
\end{example}

%Ex 7.2
\begin{example}\rm 
Consider the octahedron
\[
\Delta=\conv\set{\pm(0,1,0),\pm(2,0,0), \pm(0,0,1)}.
\]

\begin{figure}[h!]
\tdplotsetmaincoords{60}{20}
\begin{tikzpicture}[tdplot_main_coords]
\draw[thick,dotted] (-2,0,0) -- (0,1,0);
\draw[thick,dotted] (2,0,0) -- (0,1,0);
\draw[thick] (-2,0,0) -- (0,-1,0);
\draw[thick] (2,0,0) -- (0,-1,0);
%\draw[red,thick] (-1,0,0) -- (0,1,0);
%\draw[red,thick] (0,0,0) -- (0,1,0);
%\draw[red,thick] (0,1,0) -- (1,0,0);
%\draw[thick] (2,0,0) -- (-2,0,0);
\draw[thick] (0,0,1) -- (-2,0,0);
\draw[thick] (2,0,0) -- (0,0,1);
\draw[thick,dotted] (0,1,0) -- (0,0,1);
\draw[thick] (0,-1,0) -- (0,0,1);
\draw[thick] (0,0,-1) -- (-2,0,0);
\draw[thick] (2,0,0) -- (0,0,-1);
\draw[thick, dotted] (0,1,0) -- (0,0,-1);
\draw[thick] (0,-1,0) -- (0,0,-1);
%\draw[red,thick] (-1,0,0) -- (0,1,1);
%\draw[red,thick] (0,0,0) -- (0,1,1);
%\draw[red,thick] (0,1,1) -- (1,0,0);
\end{tikzpicture}
\end{figure}
%Up to isomorphism,
By symmetry, it suffices to consider two vertex cones $C(0,0,1)$, $C(2,0,0)$ by symmetry.
We construct the convex hulls $R(0,0,1)$ and $R(2,0,0)$ together with simplex triangulations inducing semi-canonical triangulations on these cones.
% for $C(0,0,1)$ (resp. for $C(2,0,0)$).

\begin{figure}[h!]
\tdplotsetmaincoords{70}{20}
\begin{tikzpicture}[tdplot_main_coords]
\draw[thick,dotted] (-2,0,0) -- (0,1,0);
\draw[thick,dotted] (2,0,0) -- (0,1,0);
\draw[thick] (-2,0,0) -- (0,-1,0);
\draw[thick] (2,0,0) -- (0,-1,0);
\draw[thick] (0,-1,0) -- (0,0,1);
\draw[red,thick] (-1,0,0) -- (0,1,0);
\draw[red,thick] (0,0,0) -- (0,1,0);
\draw[red,thick] (0,1,0) -- (1,0,0);
\draw[red, thick] (2,0,0) -- (-2,0,0);
\draw[thick] (0,0,1) -- (-2,0,0);
\draw[thick] (2,0,0) -- (0,0,1);
\draw[thick,dotted] (0,1,0) -- (0,0,1);

\draw[red,thick] (-1,0,0) -- (0,-1,0);
\draw[red,thick] (0,0,0) -- (0,-1,0);
\draw[red,thick] (0,-1,0) -- (1,0,0);
\draw[red,thick] (-1,0,0) -- (0,0,1);
\draw[red,thick] (0,0,0) -- (0,0,1);
\draw[red,thick] (0,0,1) -- (1,0,0);
%\draw[red,thick] (-1,0,0) -- (0,1,1);
%\draw[red,thick] (0,0,0) -- (0,1,1);
%\draw[red,thick] (0,1,1) -- (1,0,0);
\end{tikzpicture}
\begin{tikzpicture}[tdplot_main_coords]
%\draw[thick,dotted] (-2,0,0) -- (0,1,0);
\draw[thick,dotted] (2,0,0) -- (0,1,0);
%\draw[thick] (-2,0,0) -- (0,-1,0);
\draw[thick] (2,0,0) -- (0,-1,0);
%\draw[red,thick] (-1,0,0) -- (0,1,0);
%\draw[red,thick] (0,0,0) -- (0,1,0);
%\draw[red,thick] (0,1,0) -- (1,0,0);
%\draw[thick] (2,0,0) -- (-2,0,0);
%\draw[thick] (0,0,1) -- (-2,0,0);
\draw[thick] (2,0,0) -- (0,0,1);
\draw[thick,dotted] (0,1,0) -- (0,0,1);
\draw[thick] (0,-1,0) -- (0,0,1);
%\draw[thick] (0,0,-1) -- (-2,0,0);
\draw[thick] (2,0,0) -- (0,0,-1);
\draw[thick, dotted] (0,1,0) -- (0,0,-1);
\draw[thick] (0,-1,0) -- (0,0,-1);
\draw[thick,red] (0,-1,0) -- (1,0,0);
\draw[thick,red] (0,1,0) -- (1,0,0);
\draw[thick,red] (0,0,1) -- (1,0,0);
\draw[thick,red] (0,0,-1) -- (1,0,0);
\draw[thick,red] (2,0,0) -- (1,0,0);
%\draw[red,thick] (-1,0,0) -- (0,1,1);
%\draw[red,thick] (0,0,0) -- (0,1,1);
%\draw[red,thick] (0,1,1) -- (1,0,0);
\end{tikzpicture}
\caption{The regions $R(0,0,1)$ and $R(2,0,0)$ each consist of four simplices $S_1, \dots ,S_4$, with $(1,0,0)\in S_i$, $i=1,\dots ,4$.}
\end{figure}

It is straightforward to verify that both $C(2,0,0)$ and $C(0,0,1)$ are small cones under the induced semi-canonical triangulation. %induced by the figure. 
Moreover, the octahedron $\Delta$ can be expressed as the intersection of eight hyperplanes: $\Delta=\smallcap_{i=1}^8 H_i$ with
\begin{align*}
 H_1&=\set{(x,y,z)\in \bbR^3| x+2y+2z\leq 2}, & H_2=\set{(x,y,z)\in \bbR^3| x+2y-2z\leq 2},\\
 H_3&=\set{(x,y,z)\in \bbR^3| x-2y+2z\leq 2}, & H_4=\set{(x,y,z)\in \bbR^3| x-2y-2z\leq 2},\\
 H_5&=\set{(x,y,z)\in \bbR^3| -x+2y+2z\leq 2}, & H_6=\set{(x,y,z)\in \bbR^3| -x+2y-2z\leq 2},\\
 H_7&=\set{(x,y,z)\in \bbR^3| -x-2y+2z\leq 2}, & H_8=\set{(x,y,z)\in \bbR^3| -x-2y-2z\leq 2}.
\end{align*}
Thus, $\Delta$ is a symmetric weakly reflexive polytope. By Corollary \ref{cor:Fano},
the associated polarized toric variety $(X_\Delta, L_\Delta)$ is asymptotically Chow polystable.
\end{example}

%Ex3
\begin{example}\rm 
Consider another octahedron
\[
\Delta_3=\conv\set{\pm(0,1,0),\pm(3,0,0), \pm(0,0,1)}.
\]

\begin{figure}[h!]
\tdplotsetmaincoords{60}{20}
\begin{tikzpicture}[tdplot_main_coords]
\draw[thick,dotted] (-3,0,0) -- (0,1,0);
\draw[thick,dotted] (3,0,0) -- (0,1,0);
\draw[thick] (-3,0,0) -- (0,-1,0);
\draw[thick] (3,0,0) -- (0,-1,0);
%\draw[red,thick] (-1,0,0) -- (0,1,0);
%\draw[red,thick] (0,0,0) -- (0,1,0);
%\draw[red,thick] (0,1,0) -- (1,0,0);
%\draw[thick] (2,0,0) -- (-2,0,0);
\draw[thick] (0,0,1) -- (-3,0,0);
\draw[thick] (3,0,0) -- (0,0,1);
\draw[thick,dotted] (0,1,0) -- (0,0,1);
\draw[thick] (0,-1,0) -- (0,0,1);
\draw[thick] (0,0,-1) -- (-3,0,0);
\draw[thick] (3,0,0) -- (0,0,-1);
\draw[thick, dotted] (0,1,0) -- (0,0,-1);
\draw[thick] (0,-1,0) -- (0,0,-1);
%\draw[red,thick] (-1,0,0) -- (0,1,1);
%\draw[red,thick] (0,0,0) -- (0,1,1);
%\draw[red,thick] (0,1,1) -- (1,0,0);
\end{tikzpicture}
\end{figure}
Again, by symmetry,  %Up to isomorphism, 
it suffices to consider two cones $C(0,0,1)$ and $C(3,0,0)$.
We construct the convex hulls $R(0,0,1)$ and $R(3,0,0)$ together with the corresponding simplex triangulations
inducing semi-canonical triangulations. 
%for $C(0,0,1)$ (resp. for $C(3,0,0)$).

\begin{figure}[h!]
\tdplotsetmaincoords{70}{20}
\begin{tikzpicture}[tdplot_main_coords]
\draw[thick,dotted] (-3,0,0) -- (0,1,0);
\draw[thick,dotted] (3,0,0) -- (0,1,0);
\draw[red] (-2,0,0) -- (0,1,0);
\draw[red] (2,0,0) -- (0,1,0);
\draw[thick] (-3,0,0) -- (0,-1,0);
\draw[thick] (3,0,0) -- (0,-1,0);
\draw[red,thick] (-2,0,0) -- (0,-1,0);
\draw[red,thick] (2,0,0) -- (0,-1,0);
\draw[thick] (0,-1,0) -- (0,0,1);
\draw[red,thick] (-1,0,0) -- (0,1,0);
\draw[red,thick] (0,0,0) -- (0,1,0);
\draw[red,thick] (0,1,0) -- (1,0,0);
\draw[red, thick] (3,0,0) -- (-3,0,0);
\draw[thick] (0,0,1) -- (-3,0,0);
\draw[thick] (3,0,0) -- (0,0,1);
\draw[thick,dotted] (0,1,0) -- (0,0,1);

\draw[red,thick] (-1,0,0) -- (0,-1,0);
\draw[red,thick] (0,0,0) -- (0,-1,0);
\draw[red,thick] (0,-1,0) -- (1,0,0);
\draw[red,thick] (-1,0,0) -- (0,0,1);
\draw[red,thick] (0,0,0) -- (0,0,1);
\draw[red,thick] (0,0,1) -- (1,0,0);
%\draw[red,thick] (-1,0,0) -- (0,1,1);
%\draw[red,thick] (0,0,0) -- (0,1,1);
%\draw[red,thick] (0,1,1) -- (1,0,0);
\end{tikzpicture}
\begin{tikzpicture}[tdplot_main_coords]
%\draw[thick,dotted] (-2,0,0) -- (0,1,0);
\draw[thick,dotted] (3,0,0) -- (0,1,0);
%\draw[thick] (-2,0,0) -- (0,-1,0);
\draw[thick] (3,0,0) -- (0,-1,0);
%\draw[red,thick] (-1,0,0) -- (0,1,0);
%\draw[red,thick] (0,0,0) -- (0,1,0);
%\draw[red,thick] (0,1,0) -- (1,0,0);
%\draw[thick] (2,0,0) -- (-2,0,0);
%\draw[thick] (0,0,1) -- (-2,0,0);
\draw[thick] (3,0,0) -- (0,0,1);
\draw[thick,dotted] (0,1,0) -- (0,0,1);
\draw[thick] (0,-1,0) -- (0,0,1);
%\draw[thick] (0,0,-1) -- (-2,0,0);
\draw[thick] (3,0,0) -- (0,0,-1);
\draw[thick, dotted] (0,1,0) -- (0,0,-1);
\draw[thick] (0,-1,0) -- (0,0,-1);
\draw[thick,red] (0,-1,0) -- (2,0,0);
\draw[thick,red] (0,1,0) -- (2,0,0);
\draw[thick,red] (0,0,1) -- (2,0,0);
\draw[thick,red] (0,0,-1) -- (2,0,0);
\draw[thick,red] (3,0,0) -- (2,0,0);
%\draw[red,thick] (-1,0,0) -- (0,1,1);
%\draw[red,thick] (0,0,0) -- (0,1,1);
%\draw[red,thick] (0,1,1) -- (1,0,0);
\end{tikzpicture}
\caption{The regions $R(0,0,1)$ and $R(3,0,0)$ each consist of four simplices $S_1, \dots ,S_4$ with $(3,0,0)\in S_i$, $i=1,\dots ,4$.}\label{fig:Triangulations}
\end{figure}
Both $C(3,0,0)$ and $C(0,0,1)$ are small cones under these triangulations (see, Figure \ref{fig:Triangulations}).
Furthermore, $\Delta_3=\smallcap_{i=1}^8 H_i$, where
%the corresponding semi-canonical triangulation induced by the figure. Also, the octahedron $\Delta$ can be expressed as an intersection of eight hyperplanes: $\Delta_3=\smallcap_{i=1}^8 H_i$ with
\begin{align*}
 H_1&=\set{(x,y,z)\in \bbR^3| x+3y+3z\leq 3}, & H_2=\set{(x,y,z)\in \bbR^3| x+3y-3z\leq 3},\\
 H_3&=\set{(x,y,z)\in \bbR^3| x-3y+3z\leq 3}, & H_4=\set{(x,y,z)\in \bbR^3| x-3y-3z\leq 3},\\
 H_5&=\set{(x,y,z)\in \bbR^3| -x+3y+3z\leq 3}, & H_6=\set{(x,y,z)\in \bbR^3| -x+3y-3z\leq 3},\\
 H_7&=\set{(x,y,z)\in \bbR^3| -x-3y+3z\leq 3}, & H_8=\set{(x,y,z)\in \bbR^3| -x-3y-3z\leq 3}.
\end{align*}
Hence $\Delta_3$ is a symmetric weakly reflexive polytope. 
Therefore, by Corollary \ref{cor:Fano}, the associated polarized toric variety $(X_{\Delta_3}, L_{\Delta_3})$ is asymptotically Chow polystable.
\end{example}
Remark that \cite{LLSW19} show that $S_3=\conv\set{\pm(0,1),\pm(3,0)}$ is asymptotically Chow unstable. In \cite{Lee25}, for a polytope $\Delta$, one defines 
\[
D(\Delta):=\conv\set{(\bs p,0), \pm(0,0,1)| \bs p\in \Delta}.
\]  
Hence $D(S_3)$ is asymptotically Chow polystable, while $S_3$ itself is asymptotically Chow unstable.

%References

\end{document}